\documentclass[a4paper,11pt]{article}
\usepackage{amsmath}
\usepackage[T1]{fontenc}
\usepackage{amssymb}
\usepackage[english]{babel}
\usepackage{graphicx}

\ifx\JPicScale\undefined
\def\JPicScale{1}
\fi
\unitlength \JPicScale mm
\newtheorem{definition}{\textbf{\underline{Definition}}}
\newtheorem{lemma}{\textbf{\underline{Lemma}}}
\newtheorem{theorem}{\textbf{\underline{Theorem}}}
\newtheorem{proposition}{\textbf{\underline{Proposition}}}
\begin{document}

\begin{center}
{\Huge New existence theorems in measure theory and equivalence results for the existence of invariant probabilities}

\vspace{10 mm}

{\LARGE Larrieu Jean-Yves}

{\Large E-mail: maths.larrieu@gmail.com}

\vspace{5 mm}
\end{center}

\begin{quotation}
{\Large Abstract } We describe a construction process of a
relevant measure in any non-empty compact metric space. This probability measure has invariance properties with respect to isometric maps defined on open sets. These properties imply that this measure is an appropriate generalisation of the Lebesgue one.
Results about its uniqueness are showed, and applications and complementary properties are quickly studied. Peculiarly, we show an equivalence result in a general framework linked with the Krylov-Bogolioubov theorem.

{\Large R\'{e}sum\'{e} } Cet article dcrit un proc\'{e}d\'{e} de
construction d'une mesure int\'{e}ressante dans tout espace m%
\'{e}trique compact non vide. Elle poss\`{e}de des propri\'{e}t\'{e}%
s d'invariance vis \`{a} vis des applications isom\'{e}triques d\'{e}finies sur les ouverts. Ces propriétés montrent que cette mesure est une bonne généralisation de la mesure de Lebesgue. On prouve ensuite des théorèmes au sujet de son unicit\'{e}, puis on donne des applications de nos résultats. En particulier, on démontre un résultat d'équivalence dans un cadre proche de celui du théorème de Krylov-Bogolioubov.
\end{quotation}

\section*{Keywords:} Integral geometry, compact metric space, probability, invariant measure, measured dynamical system.

\section*{Introduction}

Presently, the most frequent point of view concerning Probability Theory is the axiomatic approach of Kolmogorov. In this framework, the space, the $\sigma$-algebra and the measure are given. Of course, some theorems show the existence of such structures,with eventually  pleasant properties. We can cite the Lebesgue measure, the Haar measures and the Hausdorff measure for example (all of them can be transformed into probability measures under good assumptions). The interesting fact about these measures is their invariance properties by some sets of maps.

Very often, situations from Physics or Geometry seem to have to be modelled by probabilities. Actually, human beings generally have intuition of probability measures that come out with the geometry of the space. In Natural Sciences, this geometry describes the physic of the situation.
The previous existence theorems seem too weak to translate this intuition in its whole generality, and concretely, one has to build the intuited probability in each case separately. For example, dynamical systems (either continuous or discrete) define sets of maps endowed with group-like structures. In such situations, one can look for probabilities invariant by this set of maps. This question naturally arises in numerous physic contemporary problems linked with Ergodic Theory, which is a very dynamical field of actual mathematical researches on his own, thanks to its fruitful applications even out of Physics. Generally, the sets of maps studied by this theory are "$1$-parameter"-sets. Taking a bigger group-like structure, for example a group whose action is transitive, showing the existence of such relevant measures can raise problems. The techniques of construction of such invariant measures are grouped together in a domain called Integral Geometry. But really few global approaches exist. This fact is a bit deceiving, and the research for global theorems appears to be relevant.

Moreover, assuming the existence of a general theorem showing those isolated cases could be treated all together, one can easily imagine that this theorem would be useful in new situations, presently out of range.

A general framework for problems of Integral Geometry is the following:

\vspace{4 mm}

{\itshape
Given a space $E$ and a groupoid $G$ of partial maps defined on $E$, is it possible to find a probability measure $m$ on a suitable $\sigma$-algebra such that $m$ is invariant by all maps of $G$?}

\vspace{4 mm}

One can observe that the three well-known previous measures fit in this frame, once forgotten the requirement that $m(E)=1$. Indeed, the Lebesgue measure on a bounded open of a real vector space answers the question for the groupoid of partial translations defined on Borel sets. The Haar measure is linked with the choice of the groupoid of left-translations (or the right ones) defined on Borel sets in a locally compact group, and the Hausdorff measure with the choice of the groupoid of partial isometries defined on open sets. We use the word "partial" to mean that the maps are not necessarily defined on the whole space. We will give precisions about this later.

In this paper, we show a theorem stating that the previous problem has solutions when $E$ is a compact metric space and $G$ is the groupoid of partial isometries defined on open sets. Some generalisations of it surely exist, but this one is interesting since it gives a new point of view about this Integral Geometry problem. It generalises the known constructions, and it has some interesting consequences. For example, we show the existence of a "nice" probability on the space of probabilities on a compact space. Moreover, one can use our theorem to show that the set of compact subsets of a compact metric space, endowed with the Hausdorff distance, can be given natural probability measures. We will quickly see that the structure of the space of these probability measures gives informations about the topological properties of the underlying compact metric space. This last application can be illustrated in the framework of the Phylogeny science. Indeed, we build a simplified model for a random choice of an evolution tree among every possible evolutions trees. This model can be used to describe how Nature chooses a distribution of mutations having as a consequence the evolution of species from a common ancestor.

\vspace{2 mm}

We also give some sufficient conditions for uniqueness of such an invariant measure, following the work of \cite{Chris}, in cases close to the situation of homogeneity of the underlying space. 

\vspace{2 mm}

Moreover, some generalisation of our main theorem are shown, for instance a Lebesgue-like measure existence and uniqueness result on locally compact homogeneous metric spaces.

Finally, a quite interesting and exhaustive theorem is given about the Krylov-Bogolioubov problem. It states results close to an equivalence for the existence of an invariant probability under the action of a monoid of measurable maps on a compact metric space. This framework is quite more general than the usual Bogolioubov theorem, since it deals with a monoid of measurable maps instead of a group of continuous ones, even if the action space is only a compact metric space. This result is another good answer to the Integral Geometry problem. 

\vspace{2 mm}

The construction process which is used in this paper to build a general probability invariant by partial isometries defined on opens is in fact very close to the original building method of Haar for his measure on locally compact groups. In both cases the axiom of choice is needed to take some kind of limits for sets functions depending on a more and more precise covering of the space. Here, we made the choice of a construction using a "limit notion". And in both cases the sets functions are defined has quotients of cardinals of these coverings. Finally, in both cases, the limit sets function is used together with suprema and infima to define a measure taking inspiration from the properties of Radon measures. Nevertheless, we choose to give a proof as constructive as possible. Indeed, the non-constructive part in entirely contained in the limit notion existence, and we build our measure using an algebra of sets, $M$, which helps, we hope, to understand how this probability works. $M$ is indeed the set of measurable sets with null boundary.

\vspace{2 mm}

As a conclusion, we give some ideas to prolong this study by the definition of new kind of probabilities, in a generalized sense.

\section{Limit notions on bounded maps space}

First of all, we define a limit notion which will be useful afterwards.

\vspace{2 mm}

\begin{definition}
Let $A$ be a sub-algebra of the set of functions on $I$, a subset of $%
\mathbf{R}$. We call limit notion in $a\in \overline{I}$ (the previous
closure has to be seen in $\overline{\mathbf{R}}$) on $A$ a morphism $L$
from the algebra $A$ to $\mathbf{R}$ such that:

- for any map $f$ of $A$, $\limsup_{x \rightarrow a} f(x) \geq L(f)\geq
\liminf_{x \rightarrow a} f(x)$,

- $L$ is increasing.
\end{definition}

\vspace{2 mm}

The interesting case for us is the following: we take for $A$ the algebra of bounded
maps from $(0;+\infty)$ to $\mathbf{R}$, and $0$ for $a$.

Let us tell some generalities about limit notions. First of all, they cannot
have any good properties with composition. Indeed, if the notion $L$
generalises the limit in $0$, for every bounded map $f$, and for every map $%
g $ with limit $0$ in $0$, we should have: $L(f\circ g)=L(f)$. Such a
situation is impossible. In order to see that, let $u(x)=\cos (1/x)$ and $%
v(x)=\sin (1/x)$. One should have $L(x\mapsto v(\frac{1}{2}x))=L(v)$, hence $%
L(2vu)=L(v)$, since $\sin (X)=2\sin (X/2)\cos (X/2)$. Then $2L(v)L(u)=L(v)$.
But $u^{2}+v^{2}=1$, hence $L(u)^{2}+L(v)^{2}=1$. Finally, $u(\frac{1}{\frac{1}{x}+%
\frac{\pi }{2}})=-v$, and $\lim_{x\rightarrow 0}\frac{1}{\frac{1}{x}+\frac{\pi }{2}}%
=0$. Therefore $L(u)=-L(v)$. Hence $|L(u)|=|L(v)|=\sqrt{2}/2$, and with $%
2L(v)L(u)=L(v)$, $|L(v)|=1$, which is absurd.

Now, we show that such limit notions exist. Unfortunately, we have to use
the axiom of choice. As a consequence, the effective calculations are quite
difficult.

In fact, we only need a limit notion in $0$ in the algebra of bounded maps from $(0;+\infty )$ to $\mathbf{R}$ for our further work. We shall prove the existence of limit notions in this special case. Our proof can be adapted to other ones.

Let $F$ be the filter of $0$ neighbourhoods in $(0;+\infty )$, and $U$
a ultra-filter containing $F$. Such a ultra-filter exists, because of Zorn Theorem. Let $f$ be a bounded map from $(0;+\infty )$ to $[-M;M]\subset 
\mathbf{R}$. Then the image filter of $U$ is a ultra-filter of $[-M;M]$, a
compact space: it converges to a real number $L(f)$. We have $%
L(f)=\lim_{U}f(x)$, with $F$ contained in $U$. $L$ is obviously an
increasing morphism of algebra. Then, if $\limsup_{x\rightarrow 0}f(x)=a$,
for any map $g$ whose limit $b$ is above $a$, there is a neighbourhood of $0$
such that $f(x)<g(x)$ on this neighbourhood. $U$ contains $F$, thus $%
L(f)\leq b$ for all $b>a=\limsup_{x\rightarrow 0}f(x)$. Hence $L(f)\leq
\limsup_{x\rightarrow 0}f(x)$. Following the same kind of proof, one shows
that $L(f)\geq \liminf_{x\rightarrow 0}f(x)$.

Hence, we have shown:

\begin{theorem}
Limit notions in $0$ in the algebra of bounded maps from $(0;+\infty )$ to $%
\mathbf{R}$ exist.
\end{theorem}

\section{Metrically compatible probabilities on compact metric spaces}

\subsection*{Notations}

Let us set a few notations. If $A$ is a subspace of a topological space, one
denote $\overline{A}$ its closure and $\overset{\circ }{A}$ its interior.
One denotes $Fr(A)$ the boundary of $A$: $Fr(A)=\overline{A}\setminus 
\overset{\circ }{A}$. In a metric space, the open ball with centre $x$ and
radius $r$ is denoted by $B(x;r)$ and the associated closed ball is $%
B_{f}(x,r)={y/d(x;y)\leq r}$. Beware to the fact that $\overline{B(x;r)}$ not necessarily equals $B_{f}(x,r)$.
We denote $\amalg$ for disjoint unions, and we say that the disjoint union $A \amalg B$ is positive if $d(A;B) > 0$. We say the same way that an inclusion $A \subset B$ is positive if $A$ and the complement of $B$ are positively disjoint.
 The complement of any subset $A$ in a set $B$ is
denoted by $\complement_{B}^A$. Finally, one denotes $\overline{P}^{\delta }$ for $%
\{x\in K/\exists y\in P/d(x;y)\leq \delta \}$, where $P$ is a subset of $K$.
We call this set the $\delta $-Minkowski sausage of $P$.

\subsection{Construction of the probability measure}

In this section, one associates a probability measure to every compact
metric space. It is such that for any couple of partially isometric opens $U$
and $V$, their measures equal. This construction only depends on
the choice of a limit notion in $0$ of the space of bounded maps from $%
(0;+\infty )$ to $\mathbf{R}$. Let us choose once for all such a limit
notion $L$. We assume that it is the same all the paper long.

First of all, one defines the notion of partially isometric opens. $\phi : U
\rightarrow V$ is a partial isometry between the opens $U$ and $V$ if it is
an isometry from $U$ to $V$. There is no need that it can be prolonged to
the whole ambient space containing $U$.

Let us now explain why the really interesting maps for invariance
properties in metric spaces are not global isometries but the partial ones,
between opens. Look at the interval $[0;1]$: there are only two global
isometries (the identity, and the symmetry around $1/2$). Nevertheless,
there are many partial isometries. Indeed, when balls are small enough, they
are partially isometric. The groupoid of partial isometries figures the
homogeneous nature of this space, except in the neighbourhood of its boundary.

Hence, the good notion is the one of partial isometries. However, if you
consider those partial isometries between any Borel sets, you can check that
there is no probability measure invariant under such maps in whole
generality. Just look at the following compact space: $\{\frac{1}{n},n\in 
\mathbf{N}\}\cup \{0\}$, with the distance associated with the absolute
value. It is denumerable, and all its singletons are partially isometric to
each other. As a consequence, the $\sigma $-additivity prevents the
existence of a suitable probability. Nevertheless, demanding only the
invariance relatively to partially isometric opens, one easily sees that
the Dirac measure in $0$ suits. One can actually show that it is the only
suitable measure in that case.

\vspace{2 mm}

Let us define a few more notations, in order to begin the construction of
our measure: $K$ is the main non-empty compact metric space, and if $A$ is
a subspace of $K$, one denotes $N(A;\varepsilon )$ the infimum of cardinals of coverings of $A$ by open balls of radius $%
\varepsilon $ and whose centres are elements of $A$. Since $K$ is compact,
this number is finite, and since it is an integer, the infimum is
reached for at least one covering. We call such a covering a $\varepsilon $%
-minimal covering of $A$, in order to sum up the exact expression which
should be "covering of minimal cardinal by open balls with radius $%
\varepsilon $ and centres chosen in $A$".

Moreover $\varepsilon \mapsto N(A;\varepsilon )$ is obviously a decreasing
map, and if $A\subseteq A^{\prime }$, $N(A;\varepsilon )\leq N(A^{\prime
};\varepsilon )$. In addition, $\varepsilon \mapsto \frac{N(A;\varepsilon )}{%
N(K;\varepsilon )}$ is bounded. Consequently, one can calculate its limit
notion $L$ in $0$.

\vspace{2 mm}

The main idea of the following construction is explained now. For any subset 
$A$ of $K$, the probability of $A$ is approximatively $\lim_{\varepsilon
\rightarrow 0}\frac{N(A;\varepsilon )}{N(K;\varepsilon )}$. Of course, the
limit notion $L$ will be useful to define this limit. In fact, one can
exhibit a compact $K$ and a subset $A$ of $K$ such that $\lim_{\varepsilon
\rightarrow 0}\frac{N(A;\varepsilon )}{N(K;\varepsilon )}$ does not exist.
Such a counter example follows. It is inspired from another one M. Drouin
showed me. The existence of such situations proves the necessity to use a
limit notion.

\vspace{2 mm}

Denote $A_0 = [0;1]$, and $A_n$ being assumed to be constructed, let us
denote $A_{n+1} = A_n \cap \lbrace x \in [0;1] / \exists p \in \mathbf{N}/
\vert x - \frac{2p+1}{13^{k+1}\times 5^k}\vert \leq \frac{1/2}{%
13^{k+1}\times 5^k} - \frac{1}{65^{k+1}}\rbrace$ if $n = 2k$ is even, and $%
A_{n+1} = A_n \cap \lbrace x \in [0;1] /\exists p \in \mathbf{N}/ \vert x - 
\frac{2p+1}{13^{k+1}\times 5^{k+1}}\vert \leq \frac{1/2}{13^{k+1}\times
5^{k+1}} - \frac{1}{65^{k+1}}\rbrace$ if $n = 2k+1$ is odd. This process
consists in two alternative elementary transformations similar to that used in the
construction of Cantor sets. One considers now $A = \bigcap_n A_n$. One shows
that $N(A;\frac{1}{13^{k+1}\times 5^k}) = 7^{k+1}\times 3^{k} $ and that $%
N(A;\frac{1}{13^k\times 5^k}) = 7^k\times 3^k$. It can be shown with the definition
of the map $N(A;.)$, taking care of the fact that for any elementary
interval of any $A_n$, its centre belongs to $A$. It follows from the choice
of the constants $13$ and $5$. (See the schema)

Now, one denotes $B_0 = [0;1]$, then $B_{2k+1} = B_{2k}\cap \lbrace x \in
[0;1] /\exists p \in \mathbf{N}/ \vert x - \frac{2p+1}{13^k\times 5^{k+1}}%
\leq \frac{1/2}{13^k\times 5^{k+1}} - \frac{1}{65^{k+1}}\rbrace$ and $B_{2k}
= B_{2k-1}\cap \lbrace x \in [0;1] / \exists p \in \mathbf{N}/ \vert x - 
\frac{2p+1}{13^{k+1}\times 5^{k+1}}\leq \frac{1/2}{13^{k+1}\times 5^{k+1}} - 
\frac{1}{65^{k+1}}\rbrace$. One has defined a decreasing sequence of
compacts with the previous process, but beginning with the other elementary
transformation. One then denotes $B = \lbrace 2\rbrace + \bigcap_n B_n$, and 
$K = A\cup B$ (it is a compact set). One shows that $N(B;\frac{1}{%
13^{k+1}\times 5^k}) = 7^k\times 3^{k+1}$ and that $N(B;\frac{1}{13^k5^k}) =
7^k\times 3^k$.

Hence $(\frac{N(A;\frac{1}{13^{k+1}\times 5^k})}{N(B;\frac{1}{13^{k+1}\times
5^k})})_k$ is a constant with value $\frac{7}{3}$, and $(\frac{N(A;\frac{1}{%
13^{k+1}\times 5^k})}{N(B;\frac{1}{13^k\times 5^k})})_k$ equals $1$.

Then, assume that $\varepsilon \mapsto \frac{N(A;\varepsilon)}{%
N(K;\varepsilon)}$ has a limit $l$ in zero. Then $N(K; \varepsilon) =
N(A;\varepsilon) + N(B;\varepsilon)$ for $\varepsilon$ small enough, because 
$A$ and $B$ are positively disjoint in $K$ (positively disjoint means that $d(A,B) > 0$). Therefore: $$\lim_{\varepsilon
\rightarrow 0} \frac{N(A;\varepsilon)}{N(K;\varepsilon)} = \lim_{\varepsilon
\rightarrow 0} \frac{N(A;\varepsilon)}{N(A;\varepsilon) + N(B;\varepsilon)}
= \lim_{\varepsilon \rightarrow 0} \frac{1}{1 + \frac{N(B;\varepsilon)}{%
N(A;\varepsilon)}} = l,$$ and $\lim_{\varepsilon \rightarrow 0} \frac{%
N(B;\varepsilon)}{N(A;\varepsilon)}$ exists, which is absurd.

\vspace{2 mm}

Consequently, it is necessary to use a generalised notion of limit in zero
in order to define $\lim_{\varepsilon \rightarrow 0} \frac{N(A;\varepsilon)}{%
N(K;\varepsilon)}$ for any subset $A$ of a general compact set $K$.

\vspace{2 mm}

We can convince us on many examples that this definition of the measure of $%
A $ is often very satisfying. Indeed, in non-empty compact metric spaces, a
probability such that many Borel sets $A$ have measure $\lim_{\varepsilon
\rightarrow 0} \frac{N(A;\varepsilon)}{N(K;\varepsilon)}$ exists. If,
unfortunately, this limit does not exist, one takes $L_{\varepsilon
\rightarrow 0} \frac{N(A;\varepsilon)}{N(K;\varepsilon)}$. Consider for
example the interval $[0;1]$ endowed with Lebesgue measure, and its
intervals.

But this limit does not suit in whole generality. To convince us, recall the
example of $K = \lbrace \frac{1}{n}, n \in \mathbf{N}^*\rbrace \cup \lbrace
0\rbrace$. As seen before, our limit does not define a measure. It follows
from the non commutation of the symbols $\lim$ and $\sum$. One has to use
this "limit" $L_{\varepsilon \rightarrow 0} \frac{N(A;\varepsilon)}{%
N(K;\varepsilon)}$ on a smaller set than $P(K)$. In fact, one will use $%
L_{\varepsilon \rightarrow 0} \frac{N(A;\varepsilon)}{N(K;\varepsilon)}$
when $A$ is open.

Before speaking of this problem, let us give a result concerning $%
N(.;\varepsilon)$. \vspace{2 mm}

\begin{lemma}
Let $A$ and $B$ be subsets of $K$. Then: $N(A;\varepsilon) +
N(B;\varepsilon) - 2N(\overline{A}^{\delta}\cap \overline{B}%
^{\delta};\varepsilon) \leq N(A\cup B; \varepsilon) \leq N(A;\varepsilon) +
N(B;\varepsilon)$ for all $\varepsilon < \delta/2$.
\end{lemma}

\vspace{2 mm}

\underline{\textbf{Demonstration:}} The second part of the inequality is
obvious. Indeed, if you merge a $\varepsilon$-minimal covering of $A$ and
one of $B$, you get a $\varepsilon$-covering of $A\cup B$.

For the other part, let us write: $A \cup B = (A \setminus \overline{B}%
^{\delta})\cup (B \setminus \overline{A}^{\delta})\cup (A \cap \overline{B}%
^{\delta})\cup (B \cap \overline{A}^{\delta}).$ Hence, with the hypothesis: 
$N(A\cup B; \varepsilon) \geq N((A \setminus \overline{B}^{\delta})\cup (B
\setminus \overline{A}^{\delta}),\varepsilon )= N(A \setminus \overline{B}%
^{\delta},\varepsilon ) + N(B \setminus \overline{A}^{\delta},\varepsilon ).$
Indeed, if a ball of a $\varepsilon$-minimal covering of $(A \setminus 
\overline{B}^{\delta})\cup (B \setminus \overline{A}^{\delta})$ has its
center in $A$, it has no intersection with $B \setminus \overline{A}%
^{\delta} $ (because $\varepsilon < \delta/2$). Hence, it is useless to
cover this second set. One can make the same remark swapping $A$ and $B$.
Then: $A = (A \setminus \overline{B}^{\delta})\cup (A \cap \overline{B}%
^{\delta})$, so $N(A;\varepsilon) \leq N(A \setminus \overline{B}%
^{\delta};\varepsilon) + N(A \cap \overline{B}^{\delta};\varepsilon)$, and,
in the same way, $N(B;\varepsilon) \leq N(B \setminus \overline{A}%
^{\delta};\varepsilon) + N(B \cap \overline{A}^{\delta};\varepsilon)$. Hence 
$N(A\cup B; \varepsilon) \geq N(A;\varepsilon) + N(B;\varepsilon) - N(A \cap 
\overline{B}^{\delta};\varepsilon) - N(B \cap \overline{A}%
^{\delta};\varepsilon) \geq N(A;\varepsilon) + N(B;\varepsilon) - 2N(%
\overline{A}^{\delta}\cap \overline{B}^{\delta})$ using the fact that $A
\mapsto N(A;\varepsilon )$ is increasing.

\begin{flushright}
$\blacksquare$
\end{flushright}

\begin{lemma}
Let $A$ and $B$ be disjoint subsets of $K$. For all $\alpha>0$, there is a $%
\delta>0$ such that $\overline{A}^{\delta}\cap \overline{B}^{\delta}
\subseteq \overline{Fr(A)}^{\alpha}$.
\end{lemma}

\vspace{2 mm}

\underline{\textbf{Demonstration:}} By absurd, if it was not, for all
natural numbers $n$ different from zero, there would be a $x_n$ in $\overline{A%
}^{1/n}\cap \overline{B}^{1/n} \setminus \overline{Fr(A)}^{\alpha}$. Then,
we could choose a subsequence of $(x_n)$ converging to some $x$ of $K$ such
that $d(x;A) = 0$, $d(x;B) = 0$ and $d(x; Fr(A)) > \alpha$ (write our
inequalities in the limit case). It is absurd.

\begin{flushright}
$\blacksquare$
\end{flushright}

\begin{definition}
One denotes $M$ the set of subsets $A$ of $K$ such that 
\begin{equation*}
\lim_{\delta \rightarrow 0} L_{\varepsilon \rightarrow 0} \frac{N(\overline{%
Fr(A)}^{\delta};\varepsilon)}{N(K;\varepsilon)}=0. 
\end{equation*}
\end{definition}

\vspace{2 mm}

\begin{proposition}
$M$ is an algebra of sets.
\end{proposition}

\vspace{2 mm}

\underline{\textbf{Demonstration:}} Since for any subset $A$, $Fr(A) =
Fr(\complement_K^A)$, $M$ is invariant by complement. Then, since for all subsets $A$
and $B$, $Fr(A\cup B) \subseteq Fr(A) \cup Fr(B)$, and for all $\delta > 0$, 
$\overline{Fr(A\cup B)}^{\delta} \subseteq \overline{Fr(A)}^{\delta} \cup 
\overline{Fr(B)}^{\delta}$, with the first lemma: 
\begin{equation*}
\frac{N(\overline{Fr(A\cup B)}^{\delta};\varepsilon)}{N(K;\varepsilon)} \leq 
\frac{N(\overline{Fr(A)}^{\delta};\varepsilon)}{N(K;\varepsilon)} + \frac{N(%
\overline{Fr(B)}^{\delta};\varepsilon)}{N(K;\varepsilon)} , 
\end{equation*}
hence 
\begin{equation*}
\lim_{\delta \rightarrow 0}L_{\varepsilon \rightarrow 0}\frac{N(\overline{%
Fr(A\cup B)}^{\delta};\varepsilon)}{N(K;\varepsilon)} 
\end{equation*}
\begin{equation*}
\leq \lim_{\delta \rightarrow 0}L_{\varepsilon \rightarrow 0} \frac{N(%
\overline{Fr(A)}^{\delta};\varepsilon)}{N(K;\varepsilon)} + \lim_{\delta
\rightarrow 0}L_{\varepsilon \rightarrow 0} \frac{N(\overline{Fr(B)}%
^{\delta};\varepsilon)}{N(K;\varepsilon)} \leq 0 , 
\end{equation*}
and one see that $M$ is stable by finite union. Moreover, it is obvious that 
$\emptyset \in M$. $M$ is an algebra of sets.

\begin{flushright}
$\blacksquare$
\end{flushright}

$A \mapsto
L_{\varepsilon\rightarrow 0}(\frac{N(A;\varepsilon)}{N(K;\varepsilon)})$,
defined on $M$, is nearly a measure. Indeed, denoting this map $m$:

\begin{proposition}
$m$ is finitely additive in $M$.
\end{proposition}

\vspace{2 mm}

\underline{\textbf{Demonstration:}} Let $A$ and $B$ be disjoint sets
belonging to $M$. Using the lemma $2$, for all $\alpha>0$, there is a $%
\delta>0$ such that $\overline{A}^{\delta}\cap \overline{B}^{\delta}
\subseteq \overline{Fr(A)}^{\alpha}$.

Then: 
\begin{equation*}
N(A;\varepsilon) + N(B;\varepsilon) - 2N(\overline{Fr(A)}^{\alpha};%
\varepsilon) \leq N(A\cup B; \varepsilon) \leq N(A;\varepsilon) +
N(B;\varepsilon) 
\end{equation*}
for all $\alpha>0$ and all $\varepsilon >0$ small enough. Then, you just
have to divide by $N(K;\varepsilon)$, to take $L_{\varepsilon \rightarrow 0}$%
, and then to use $\lim_{\alpha \rightarrow 0}$. So we have $m(A \cup B) =
m(A) + m(B)$.

\begin{flushright}
$\blacksquare$
\end{flushright}

But, at this level, we are still far away from having built a real measure.
Indeed, the $\sigma$-additivity in $M$ seems to be difficult to prove.
Moreover, the following lemma insures us that the behaviour of this map of
sets $m$ does not suit.

\begin{lemma}
Let $P$ be in $M$. Then $\overset{\circ}{P}$ and $\overline{P}$ are in $M$
and $m(P) = m(\overset{\circ}{P}) = m(\overline{P}).$
\end{lemma}

\vspace{2 mm}

\underline{\textbf{Demonstration:}} If $P$ is in $M$, $Fr(P)$ is in $M$,
hence $\overset{\circ}{P}$ and $\overline{P}$ too, since $M$ is an algebra
of sets. Then, since $m(Fr(P)) =0$, $m(\overline{P}) = m(\overset{\circ}{P})
\leq m(P) \leq m(\overline{P})$, which leads to the conclusion.

\begin{flushright}
$\blacksquare$
\end{flushright}

This lemma, and the one in the same spirit which will follow, will be useful
in the future.

\vspace{2 mm}

In order to see that this result prevents $m$ from having any extension to
the expected measure, let us consider one more time $K = \lbrace \frac{1}{n}, n \in 
\mathbf{N}\rbrace \cup \lbrace 0\rbrace$, and $P = \lbrace 0\rbrace$. The
measure that we are looking for (the Dirac in $0$) does not satisfy the
condition of the lemma.

Moreover, there is another issue linked with the definition of the algebra
of sets $M$: the measure should have properties using the opens of $K$, and
one can easily make situations where some opens of $K$ are not in $M$. For
example, take the interval $[0;1]$. The expected measure is the Lebesgue
measure on $[0;1]$ (you can show it by cutting $[0;1]$ into several
intervals with dyadic ends, and calculating the measure of the union of
disjoint intervals). Observe the open set $O$ taken to be the complement of a
non-zero Lebesgue measure Cantor set. $O$ has then a boundary of
non-zero measure, and as a consequence, it does not belong to $M$. One
solves this issue with the following method:

\begin{definition}
One defines for any subset $P$ of $K$: 
\begin{equation*}
m^{\prime }(P) = \sup_{U} m(U), 
\end{equation*}
where the supremum is taken for the elements $U$ of $M$ such that $%
U\subseteq P$.
\end{definition}

\vspace{2 mm}

$m^{\prime }$ is clearly increasing.

\vspace{2 mm}

\begin{lemma}
$m^{\prime }(K) = 1$.
\end{lemma}

\vspace{2 mm}

\underline{\textbf{Demonstration:}} $m^{\prime }(K)\geq 1$
since $K$ is open. Then $m^{\prime }(K)\leq 1$ since it is the supremum of numbers smaller than $1$%
.

\begin{flushright}
$\blacksquare$
\end{flushright}

Here is a result showing that we are working in the good direction.

\begin{proposition}
$m^{\prime }$ is additive on finite families of disjoint open sets of $K$.
\end{proposition}

\vspace{2 mm}

\underline{\textbf{Demonstration:}} Let $O$ and $O^{\prime }$ be disjoint
opens of $K$. Let $U$ and $U^{\prime }$ be elements of $M$ such that $%
U\subseteq O$ and $U^{\prime }\subseteq O^{\prime }$, with $m^{\prime }(O)
\leq m(U) +\varepsilon$ and $m^{\prime }(O^{\prime }) \leq m(U^{\prime })
+\varepsilon$ for a chosen $\varepsilon>0$. Then $U\cup U^{\prime }$ is an
element of $M$ included in $O\cup O^{\prime }$. Hence: $m^{\prime }(O\cup
O^{\prime }) \geq m(U) + m(U^{\prime }) \geq m^{\prime }(O) + m^{\prime
}(O^{\prime }) - 2\varepsilon.$ Since $\varepsilon$ can be freely chosen, we
have a part of the equality.

Now, we show the sub-additivity of $m^{\prime }$ for disjoint opens of $K$.
Let $U$ be an element of $M$ such that $U \subseteq O\cup O^{\prime }$ and $%
m^{\prime }(O\cup O^{\prime }) \leq m(U) + \varepsilon$. The following proposition
shows that $O\cap U$ and $O^{\prime }\cap U$ belong to $M$. Hence $m^{\prime
}(O\cup O^{\prime }) \leq m(U) + \varepsilon \leq m((O\cap U)\amalg
(O^{\prime }\cap U)) + \varepsilon \leq m(O\cap U) + m(O^{\prime }\cap U) +
\varepsilon \leq m^{\prime }(O) + m^{\prime }(O^{\prime }) + \varepsilon.$
Note that we use the symbol $\amalg$ for disjoint unions. Since any $%
\varepsilon>0$ can be chosen, we deduct the expected result.

\begin{flushright}
$\blacksquare$
\end{flushright}

The following proposition takes a central place in the construction of our measure.

\vspace{2 mm}

\begin{proposition}
Let $O$ and $O^{\prime }$ be disjoint opens of $K$ and $P$ an element of $M$
included in $O\amalg O^{\prime }$. Then $O\cap P$ belongs to $M$ (and $%
O^{\prime }\cap P$ too).
\end{proposition}

\vspace{2 mm}

\underline{\textbf{Demonstration:}} Under those hypothesis, one has $%
Fr(O\cap P) \subseteq Fr(P)$. Indeed, choose a $x$ in $Fr(O\cap P)$. Then $x
\in \overline{O\cap P} \subseteq \overline{P}$. Show now that if $x \in 
\overset{\circ}{P}$ and $x \in \overline{O\cap P}$, $x \in \overbrace{P\cap O%
}^{\circ} = \overset{\circ}{P} \cap \overset{\circ}{O} = \overset{\circ}{P}
\cap O$. We just have to show that $x \in O$. $x \in \overline{O\cap P}\cap 
\overset{\circ}{P} \subset \overline{O}\cap \overset{\circ}{P} \subset 
\overline{O}\cap (O \cup O^{\prime }) = (\overline{O} \cap O) \cup (%
\overline{O} \cap O^{\prime }) = O\cup (\overline{O} \cap O^{\prime })$. But 
$O$ and $O^{\prime }$ are disjoint opens, therefore $\overline{O} \cap
O^{\prime }= \emptyset$.

Hence $Fr(O\cap P) \subseteq Fr(P)$. Since $P$ is in $M$, $O\cap P$ is in $M$
too.

\begin{flushright}
$\blacksquare$
\end{flushright}

It seems that our construction has a better behaviour with open sets.
Nevertheless, the following lemma shows that the aim is not reached, since $%
m^{\prime }$ is still not the good measure.

\begin{lemma}
Let $P$ be a subset of $K$. Then $m^{\prime }(P) = m^{\prime }(\overset{\circ%
}{P}) = m^{\prime }(\overline{P}).$
\end{lemma}

\vspace{2 mm}

\underline{\textbf{Demonstration:}} First of all, since $m^{\prime }$ is
increasing: $m^{\prime }(\overset{\circ}{P}) \leq m^{\prime }(P) \leq
m^{\prime }(\overline{P}).$ Choose $\varepsilon>0$, and $U$ an element of $M$
such that $U \subseteq \overline{P}$ and $m^{\prime }(\overline{P}) \leq
m(U) + \varepsilon$. With the following lemma, one can see that $%
Fr(U\setminus Fr(P)) \subseteq Fr(U)$. Therefore $U \setminus Fr(P)$ is in $%
M $ and $m(U \setminus Fr(P)) = m(U) - m(U \cap Fr(P)) = m(U) $ since $m(U
\cap Fr(P)) \leq m(Fr(P)) =0$. Hence $m^{\prime }(%
\overline{P}) \leq m(U \setminus Fr(P)) + \varepsilon$. But $U \setminus
Fr(P)$ is an element of $M$ included in $\overset{\circ}{P}$. Hence: $%
m^{\prime }(\overline{P}) \leq m^{\prime }(\overset{\circ}{P}) + \varepsilon$%
. Since $\varepsilon>0$ can be freely chosen, $m^{\prime }(\overline{P})
\leq m^{\prime }(\overset{\circ}{P})$, which implies the expected conclusion.

\begin{flushright}
$\blacksquare$
\end{flushright}

\begin{lemma}
Let $P$ be a subset of $Q$. Then $P\cap Fr(Q) \subseteq Fr(P)$.
\end{lemma}

\vspace{2 mm}

\underline{\textbf{Demonstration:}} First of all, it is obvious that $P\cap
Fr(Q) \subseteq \overline{P}$. Then, let $x$ be in $\overset{\circ}{P}$.
Then $x \in \overset{\circ}{Q}$, hence $x\not\in Fr(Q)$, and therefore $%
x\not\in P\cap Fr(Q)$.

\begin{flushright}
$\blacksquare$
\end{flushright}

Actually, one has another interesting result about the additivity on opens.
It is a very important property for the end of the construction. The
demonstration uses the fact that for opens of $K$, $m^{\prime }(O)=m^{\prime
}(\overline{O})$.

\vspace{2 mm}

\begin{lemma}
Let $O$ and $O^{\prime }$ be opens of $K$. Then $m^{\prime }(O\cup O^{\prime
}) \leq m^{\prime }(O) + m^{\prime }(O^{\prime })$.
\end{lemma}

\vspace{2 mm}

\underline{\textbf{Demonstration:}} 
\begin{equation*}
m^{\prime }(O)+m^{\prime }(O^{\prime }) \geq 
\end{equation*}
\begin{equation*}
m^{\prime }(O)+m^{\prime }(\overbrace{O^{\prime }\setminus O}^{\circ}) =
m^{\prime }(O \amalg \overbrace{O^{\prime }\setminus O}^{\circ}) = m^{\prime
}(\overline{O \amalg \overbrace{O^{\prime }\setminus O}^{\circ}})\geq 
\end{equation*}
\begin{equation*}
m^{\prime }(\overline{O\cup O^{\prime }}) = m^{\prime }(O\cup O^{\prime }). 
\end{equation*}
Indeed, $O$ and $\overbrace{O^{\prime }\setminus O}^{\circ}$ are disjoint
opens, and $\overline{O\cup O^{\prime }} \subseteq \overline{O \amalg 
\overbrace{O^{\prime }\setminus O}^{\circ}}$.

Let us show this last affirmation. First, obviously: $\overline{O} \subseteq \overline{O \amalg 
\overbrace{O^{\prime }\setminus O}^{\circ}}$. Let $x$ be in $\overline{%
O^{\prime }}\setminus \overline{O}$. There is a ball of centre $x$ and
positive radius included in the complement of $O$, and a sequence $(x_n)$ of
points of $O^{\prime }$ converging to $x$. Hence, from a certain rank, the $%
x_n$ belong to the previous ball. Therefore $x \in \overline{O^{\prime
}\setminus \overline{O}} = \overline{\overbrace{O^{\prime }\setminus O}%
^{\circ}}$.

\begin{flushright}
$\blacksquare$
\end{flushright}

Now, we reach the last part of our construction. The map of sets $m^{\prime }
$ seems to have good properties with respect to opens of $K$, but it does
not suit well with the other sets. But the expected measure is a probability
in a metric space defined on Borel sets. Therefore, it is regular, which
implies that knowing it on opens enables us to know it on the whole Borel $%
\sigma$-algebra. Following the construction of the Hausdorff measure, we put the
definition:

\begin{definition}
For any subset $P$ of $K$: 
\begin{equation*}
\mu(P) = \inf_{E_i} \sum_{i} m^{\prime }(E_i), 
\end{equation*}
where the infimum is taken on the countable families of opens $((E_i)_i)$%
such that $P\subseteq\bigcup_i E_i$.
\end{definition}

\vspace{2 mm}

Recall that it means that:

\begin{equation*}
\mu(P) = \inf_{E_i} \sum_{i} \sup_{U \in M, U \subset E_i} L_{\varepsilon
\rightarrow 0}\frac{N(U;\varepsilon)}{N(K;\varepsilon)}. 
\end{equation*}

\vspace{2 mm}

$\mu$ is an exterior measure (see \cite[page 194]{Fractales}):

\vspace{2 mm}

\begin{definition}
If $E$ is a metric space, an exterior measure $\mu$ on $E$ is a map from $%
P(E)$, the set of all subsets of $E$, in $[0;+\infty]$ such that:

- $\mu(\emptyset) =0$,

- $\mu(\bigcup_n E_n) \leq \sum_n \mu (E_n)$ for any countable family $E_n$
of subsets of $E$,

- $\mu$ is increasing.
\end{definition}

\vspace{2 mm}

\begin{proposition}
If $F$ is a covering of a metric space $E$, and if $f: F \rightarrow
[0;+\infty]$ is a map, then the map $\mu$ defined by:

- $\mu (A) = \inf \lbrace \sum_n f(F_n), F_n \in F, A \subseteq \bigcup_n
F_n\rbrace,$ if $A \neq \emptyset$, where $(F_n)$ is countable,

- $\mu (\emptyset) = 0$,

is an exterior measure.
\end{proposition}

\vspace{2 mm}

We will then use the metric exterior measures theorem which can be read in 
\cite[page 195]{Fractales}. In order to do this, we have to show a few lemmas.

\begin{lemma}
$\mu(K) = 1.$
\end{lemma}

\vspace{2 mm}

\underline{\textbf{Demonstration:}} First of all, it is obvious that $%
\mu(K) \leq m^{\prime }(K) =1$, since $K$ is an open containing $K$. Choose
a $\varepsilon>0$. If $(E_i)_i$ is a countable open covering of $K$ such
that $\mu(K) \geq \sum_{i} m^{\prime }(E_i) - \varepsilon$, since $K$ is
compact, one can extract a finite sub-covering $(O_i)$ from $(E_i)_i$. Then $%
\mu(K)+\varepsilon \geq\sum_{i} m^{\prime }(E_i) \geq \sum_{i} m^{\prime
}(O_i) \geq m^{\prime }(\bigcup_{i} O_i) = m^{\prime }(K) =1$, using the
lemmas $4$ and $6$. Hence, since $\varepsilon$ is freely chosen, $\mu(K)\geq
1$.

\begin{flushright}
$\blacksquare$
\end{flushright}

\begin{lemma}
$\mu$ is metric, which means that if $P$ and $P^{\prime }$ are positively
disjoint subsets, $\mu(P\cup P^{\prime }) = \mu(P) + \mu(P^{\prime })$.
\end{lemma}

\vspace{2 mm}

\underline{\textbf{Demonstration:}} Since $\mu$ is an exterior measure, we
just have to show that $\mu(P\cup P^{\prime }) \geq \mu(P) + \mu(P^{\prime
}).$ Choose a $\varepsilon>0$, and $(E_i)_i$ a countable open covering of $%
P\cup P^{\prime }$ such that $\mu(P\cup P^{\prime }) \geq \sum_{i} m^{\prime
}(E_i) - \varepsilon$. $P$ and $P^{\prime }$ are positively disjoint, hence
there are two disjoint opens $O$ and $O^{\prime }$ containing respectively $%
P$ and $P^{\prime }$. Then, for all $i$, $E_i \cap O$ and $E_i \cap
O^{\prime }$ are open. $P$ is covered by $(E_i \cap O)_i$, and $P^{\prime }$
is covered by $(E_i \cap O^{\prime })_i$. Therefore: $\mu(P\cup P^{\prime
}) \geq \sum_{i} m^{\prime }(E_i) - \varepsilon \geq \sum_{i} m^{\prime
}((E_i \cap O)\cup (E_i \cap O^{\prime })) - \varepsilon = \sum_{i}
m^{\prime }(E_i \cap O) + \sum_{i} m^{\prime }(E_i \cap O^{\prime }) -
\varepsilon \geq \mu(P) + \mu(P^{\prime }) - \varepsilon$. Since it is true
for all $\varepsilon>0$, we get the expected result.

\begin{flushright}
$\blacksquare$
\end{flushright}

Thus, $\mu$ is an exterior measure on $K$. Using the exterior measures
theorem, we can state that its restriction to Borel sets is a measure.
Recall here this theorem shown in \cite{Fractales}[pages 195-199].

\vspace{2 mm}

\begin{theorem}
Let $E$ be a metric space. If $\mu$ is a metric exterior measure on $E$, its
restriction to the $\sigma$-algebra of Borel sets is a measure.
\end{theorem}

\vspace{2 mm}

Moreover, $\mu(K) = 1$, hence $\mu$ is a probability on $K$.

\vspace{2 mm}

\begin{theorem}
On every non-empty metric compact space, there is a probability defined on
Borel sets such that for any couple of partially isometric open sets $U$ and 
$V$, the probability of $U$ and $V$ equal. One denotes it $\mu_K^L$, or $%
\mu_K$ if the choice of the limit notion $L$ is fixed once for all. This
construction process is canonical, once the limit notion $L$ is chosen.
Indeed, if $K$ and $K^{\prime }$ are isometric non-empty metric compact
spaces, the two measures built on $K$ and $K^{\prime }$ are image of each
other by that isometry. Moreover, if $U$ and $V$ are isometric opens of $K$,
and $A$ and $B$ two Borel sets of respectively $U$ and $V$ images of each
other by that isometry, their measures equal.
\end{theorem}

\vspace{2 mm}

\underline{\textbf{Demonstration:}} We just have to show the invariance
property with respect to isometries. Let $\phi$ be an isometry between two
compacts $K$ and $K^{\prime }$. It is obvious that if $A$ and $B$ are
partially isometric subsets of respectively $K$ and $K^{\prime }$ via $\phi$%
, for all $\varepsilon>0$, $N(A;\varepsilon) = N(B;\varepsilon)$ (use a
double inequality). Hence $M_K$ and $M_{K^{\prime }}$ are bijectively in
relation by $\phi$, and if $A$ is in $M_K$, $m_K(A) = m_{K^{\prime
}}(\phi(A))$. Then, with a double inequality, $m^{\prime }_K(O) = m^{\prime
}_{K^{\prime }}(\phi(O))$ for all opens $O$ of $K$. Then $\mu_K(F) =
\mu_{K^{\prime }}(\phi(F))$ for all subsets $F$ of $K$ (a double inequality
again), so $\mu_{K^{\prime }}$ and the push-forward measure of $\mu_{K}$
associated with $\phi^{-1}$ equal.

Moreover, if $U$ and $V$ are partially isometric opens of $K$ via $\psi$,
we can show the same way that if $\varepsilon< d(F;\complement_K^U)$, 
\begin{equation*}
N(F;\varepsilon) = N(\psi (F);\varepsilon), 
\end{equation*}
where $F$ is a subset positively contained in $U$. "positively contained"
means here that $d(F;\complement_K^U)>0$. The $\varepsilon$-coverings $F$ and $\psi (F)$
are indeed isometric for $\varepsilon$ small enough.

Hence, with $A$ positively contained in $U$, $A\in M_K \Leftrightarrow \psi
(A) \in M_K$ (the boundaries $Fr(A)$ et $Fr(\psi (A))$ are positively
contained in respectively $U$ and $V$) and $m_K(A) = m_K(\psi (A))$.

With a double inequality, one shows that if $O$ is positively contained in $%
U$, $m^{\prime }(O) = m^{\prime }(\psi (O))$.

Then, for all open sets $O$ positively contained in $U$, 
\begin{equation*}
\mu(O) = \inf_{E_i} \sum_{i} m^{\prime }(E_i) \geq \inf_{E_i} \sum_{i}
m^{\prime }(E_i \cap O), 
\end{equation*}
where $(E_i)$ is a countable family of opens such that $%
O\subseteq\bigcup_i E_i$. But then, $(E_i \cap O)$ is a countable family
of opens such that $O=\bigcup_i E_i \cap O$. Therefore, the calculation of $%
\mu(O)$ for any open $O$ can be made with the $E_i$ contained in $O$.
Therefore $\mu(O) = \mu(\psi(O))$ with a double inequality for the open sets 
$O$ positively contained in $U$.

Let us now show that $U$ and $V$ have the same $\mu$-measure. It is enough
to show that for all $\alpha>0$, there is a $r>0$ such that $\mu
(P_r)<\alpha $ and $\mu (\psi (P_r))<\alpha$, where $P_r=\lbrace x \in K /
0<d(x;\complement_K^U)\leq r\rbrace$. Indeed, $P_r$ is measurable as an intersection
of an open with a closed, and $U\setminus P_r$ is an open
positively contained in $U$.

$\bigcap_{n} P_{1/n}$ and $\bigcap_{n} \psi (P_{1/n})$ are decreasing empty
intersections of Borel sets, hence the sequences $\mu(P_{1/n})$ and $%
\mu(\psi (P_{1/n}))$ converge to $0$. It implies that $\mu(U) = \mu(V)$ (recall that $V = \psi (U)$, with $\psi$ a partial isometry defined on an open $U$ of $K$).

Now, let us finish by showing that under the theorem assumptions the Borel subsets 
$A$ and $B$ have the same measure. Since $\mu_K$ is a probability measure
defined on Borel sets, it is regular, hence $\mu_K (A) = \inf_{A \subseteq O open} \mu_K (O),$ and $\mu_K (B) = \inf_{B \subseteq O^{\prime } open} \mu_K (O^{\prime })$. Considering opens under the form $O\cap
U $ and $O^{\prime }\cap V$, one can write:

\begin{center}
$\mu_K (A) = \inf_{A \subseteq O \subseteq U} \mu_K (O)$, and $\mu_K (B) =
\inf_{B \subseteq O^{\prime }\subseteq V} \mu_K (O^{\prime })$.
\end{center}

Then, using the previous invariance property for isometric opens: 
\begin{equation*}
\mu_K (A) = \inf_{A \subseteq O \subseteq U} \mu_K (O) = \inf_{A \subseteq O
\subseteq U} \mu_K (\psi (O)) \geq \inf_{B \subseteq O^{\prime }\subseteq V}
\mu_K (O^{\prime }) = \mu_K (B), 
\end{equation*}
where $\psi$ is the partial isometry between $U$ and $V$. Reversing the roles played by $A$ and $B$, we get the inverse inequality, which implies the conclusion.

\begin{flushright}
$\blacksquare$
\end{flushright}

\subsection{Some properties satisfied by our probabilities}

Let us give now some properties of our measures. More exactly, we give some
properties of the functor of the category of non-empty compact metric spaces
to the category of non-empty compact metric spaces endowed with the
probability we built, a limit notion $L$ being fixed. We will precise later
the good morphism notion for both categories.

\subsubsection{Some inequalities for Borel sets. Some consequences for the
algebra $M$.}

\vspace{2 mm}

Our construction quite satisfies us, but it has an important drawback: in
concrete situations, when a particular space $K$ is chosen, the calculation
of the measure seems difficult. Indeed, the use of the definition is
generally too complex to be effective. We give now some elements to make
this calculation easier.

\vspace{2 mm}

\begin{proposition}
Using the notations of the beginning of this part, we have for all closed
sets $A$ of $K$: 
\begin{equation*}
L_{\varepsilon\rightarrow 0} \frac{N(A;\varepsilon)}{N(K;\varepsilon)} \leq
\mu_K(A). 
\end{equation*}
Moreover, for all opens $O$ of $K$, 
\begin{equation*}
L_{\varepsilon\rightarrow 0} \frac{N(O;\varepsilon)}{N(K;\varepsilon)}\geq
\mu_K(O). 
\end{equation*}
\end{proposition}

\vspace{2 mm}

\underline{\textbf{Demonstration:}} Let us first show that for all opens $O$, 
\begin{equation*}
L_{\varepsilon\rightarrow 0} \frac{N(O;\varepsilon)}{N(K;\varepsilon)}\geq
\mu_K(O). 
\end{equation*}
It is enough to choose in the definition of $\mu_K(O)$ the covering
containing $O$ alone. Hence $\mu_K(O) \leq m^{\prime }(O)$. Then, it is
obvious that $m(U) \leq L_{\varepsilon\rightarrow 0} \frac{N(O;\varepsilon)}{%
N(K;\varepsilon)}$ for all $U$ in $M$ contained in $O$, hence 
\begin{equation*}
m^{\prime }(O) \leq L_{\varepsilon\rightarrow 0} \frac{N(O;\varepsilon)}{%
N(K;\varepsilon)}, 
\end{equation*}
taking the supremum. Therefore $L_{\varepsilon\rightarrow 0} \frac{%
N(O;\varepsilon)}{N(K; \varepsilon)}\geq \mu_K(O).$

Let us show now the other inequality. Let $A$ be a closed. One has 
\begin{equation*}
\bigcap_{n \in \mathbf{N}} \overline{A}^{1/n} = A, 
\end{equation*}
and all sets in this formula are closed. Thus, for all $\delta>0$,
there is a $\alpha>0$ such that $\mu_K(\overline{A}^\alpha) \leq \mu_K(A) +
\delta$. And $\complement_K^{\overline{A}^\alpha}$ is open. Hence: 
\begin{equation*}
L_{\varepsilon\rightarrow 0}\frac{N(\complement_K^{\overline{A}^\alpha};\varepsilon)}{%
N(K;\varepsilon)}\geq \mu_K(\complement_K^{\overline{A}^\alpha}) = 1 - \mu_K(\overline{%
A}^\alpha) \geq 1-\delta-\mu_K(A). 
\end{equation*}

But, if $\varepsilon < \alpha/2$, $N(\complement_K^{\overline{A}^\alpha};\varepsilon)
+N(A;\varepsilon) \leq N(K;\varepsilon)$, hence: 
\begin{equation*}
\frac{N(\complement_K^{\overline{A}^\alpha};\varepsilon)}{N(K;\varepsilon)}\leq 1 - 
\frac{N(A;\varepsilon)}{N(K;\varepsilon)}. 
\end{equation*}
Taking the limit notion as $\varepsilon$ approaches $0$, 
\begin{equation*}
1-\delta-\mu_K(A) \leq L_{\varepsilon\rightarrow 0}\frac{N(\complement_K^{\overline{A}%
^\alpha};\varepsilon)}{N(K;\varepsilon)}\leq 1 - L_{\varepsilon\rightarrow
0} \frac{N(A;\varepsilon)}{N(K;\varepsilon)}. 
\end{equation*}
Hence: $L_{\varepsilon\rightarrow 0} \frac{N(A;\varepsilon)}{%
N(K;\varepsilon)} \leq \delta+\mu_K(A)$, and since it is true for all $%
\delta > 0$, we get the expected inequality.

\begin{flushright}
$\blacksquare$
\end{flushright}

\begin{proposition}
Still using the notations of this part, for all subsets $P$ in $K$: 
\begin{equation*}
P \in M \Leftrightarrow \mu_K(Fr(P)) = 0. 
\end{equation*}
\end{proposition}

\vspace{2 mm}

\underline{\textbf{Demonstration:}} Let us show the direction $\Rightarrow$: let $P$ be in $M$. Then 
\begin{equation*}
\lim_{\delta \rightarrow 0} L_{\varepsilon \rightarrow 0} \frac{N(\overset{%
\circ}{\overline{Fr(P)}^\delta};\varepsilon)}{N(K;\varepsilon)}=0. 
\end{equation*}
But $\overset{\circ}{\overline{Fr(P)}^\delta}$ is an open containing $Fr(P)$%
, hence with the previous property, 
\begin{equation*}
\mu_K(Fr(P)) \leq \mu_K(\overset{\circ}{\overline{Fr(P)}^\delta}) \leq
L_{\varepsilon \rightarrow 0} \frac{N(\overset{\circ}{\overline{Fr(P)}^\delta%
};\varepsilon)}{N(K;\varepsilon)} 
\end{equation*}
for all $\delta>0$, which implies the conclusion when $\delta$ approaches
zero.

Now, let us show the other direction: let $P$ be a subset of $K$ such that $%
\mu_K(Fr(P)) = 0.$ One has $\bigcap_{\delta >0} \overline{Fr(P)}^{\delta} =
Fr(P)$. Hence for all $\alpha>0$, there is a $\delta >0$ such that $\mu_K(%
\overline{Fr(P)}^{\delta}) < \alpha$. Now $\overline{Fr(P)}^{\delta}$ is
closed, hence with the previous property, $\alpha > \mu_K(\overline{Fr(P)}%
^{\delta}) \geq L_{\varepsilon \rightarrow 0} \frac{N(\overline{Fr(P)}%
^{\delta};\varepsilon)}{N(K;\varepsilon)}$. Therefore 
\begin{equation*}
\lim_{\delta \rightarrow 0} L_{\varepsilon \rightarrow 0} \frac{N(\overline{%
Fr(P)}^\delta; \varepsilon)}{N(K;\varepsilon)}=0, 
\end{equation*}
and $P$ belongs to $M$.

\begin{flushright}
$\blacksquare$
\end{flushright}

Hence, we can find elements of $M$ with the following method: choose a
candidate, say, $A$. Then, calculate its boundary, cover it with a countable
family of pieces, such as each one can be plunged via partial isometries a
number as tall as wished in $K$, in a disjoint way. Since $\mu_K$ exists, it
shows that each piece is a null set, hence that the boundary of $A$ is null
too. It means that $Fr(A)$ is an intrinsically null set (see later). It
follows that this subset of $K$ is in $M$. Actually, the calculation of $%
\mu_K(P)$ is really simpler when $P \in M$.

\begin{proposition}
With the notations of this paragraph, if $P$ belongs to $M$, $\mu_K(P) =
L_{\varepsilon\rightarrow 0} \frac{N(P;\varepsilon)}{N(K;\varepsilon)}$.
\end{proposition}

\vspace{2 mm}

\underline{\textbf{Demonstration:}} If $P$ belongs to $M$, $\mu_K(Fr(P)) =
0 $. Hence: 
\begin{equation*}
\mu_K(P) = \mu_K(\overset{\circ}{P}) \leq L_{\varepsilon\rightarrow 0} \frac{%
N(\overset{\circ}{P};\varepsilon)}{N(K;\varepsilon)}\leq
L_{\varepsilon\rightarrow 0} \frac{N(P;\varepsilon)}{N(K;\varepsilon)}, 
\end{equation*}
and 
\begin{equation*}
\mu_K (P)=\mu_K(\overline{P}) \geq L_{\varepsilon\rightarrow 0} \frac{N(%
\overline{P};\varepsilon)}{N(K;\varepsilon)}\geq L_{\varepsilon\rightarrow
0} \frac{N(P;\varepsilon)}{N(K;\varepsilon)}. 
\end{equation*}
It implies the expected equality.

\begin{flushright}
$\blacksquare$
\end{flushright}

Now, we have to show that there are many enough elements in $M$.

\begin{lemma}[Density]
Let $K$ be a non-empty compact metric space, and $P$ a subset of $K$. Then $%
\lbrace r>0/ \lbrace x \in P/ d(x;\complement_K^P) \geq r\rbrace \in M_K\rbrace$ is
the complement of a countable set. In the same way, $\lbrace r>0/ \overline{P%
}^r \in M_K\rbrace$ is the complement of a countable set. Consequently, for
all inclusions $F \subseteq O$ where $F$ is closed and $O$ open, there is an
element $P$ of $M_K$ contained in $O$ and containing $F$, and one can
choose $P$ open or closed.
\end{lemma}

\vspace{2 mm}

\underline{\textbf{Demonstration:}} Denote $F_r = \lbrace x \in P/
d(x;\complement_K^P) = r\rbrace$. Put, for all $\varepsilon>0$, $E_\varepsilon =
\lbrace r >0/ \mu(F_r)>\varepsilon\rbrace$.

$E_\varepsilon$ is finite. Indeed, if $E_\varepsilon$ contains more than $%
E(1/\varepsilon)+1$ elements, since the $F_r$ are disjoint, $\mu(K) = 1 \geq
\sum_{r \in E_\varepsilon} \mu(F_r) \geq (E(1/\varepsilon)+1)\varepsilon > 1$%
, which is absurd.

Hence, since for all $r$, $Fr(\lbrace x \in P/ d(x;\complement_K^P) \geq r\rbrace)
\subseteq F_r$, we conclude that the elements of $\lbrace r>0/ \lbrace x \in
P/ d(x;\complement_K^P) \geq r\rbrace \in M_K\rbrace$ are in $\complement_{(0;+\infty)}^{%
\bigcup_n E_{1/n}}$ which is the complement of a countable set.

For the last property, if $O \neq K$, one put $\delta = d(F, \complement_K^O)>0$. It
is enough to choose a $r$ in $\complement_{(0;+\infty)}^{\bigcup_n E_{1/n}} \cap
(0;\delta)$ and to take $P = \lbrace x \in O/ d(x;\complement_K^O) \geq r\rbrace$ if you wish a closed, or $%
P = \lbrace x \in O/ d(x;\complement_K^O) > r\rbrace$ if you wish an open set. ($x \mapsto d(x;\complement_K^O)$ is indeed definite and continuous
as soon as $\complement_K^O \neq \emptyset)$. If $O = K$, one considers the case when $%
F = K$ (no problem since $K \in M$), and the case $F \neq K$. In the last
one, just choose a $x \in K \setminus F$ and use the previous result with $F$
and the open $K \setminus \lbrace x\rbrace$.

Let us now show the result about $\lbrace r>0/ \overline{P}^r \in M_K\rbrace$%
. $K \setminus \overline{P}$ is open, and for all $r>0$, $\complement_K^{\overline{P}%
^r} = \complement_K^{\lbrace x \in K/ d(x;P) \leq r\rbrace} = \lbrace x \in \complement_K^P/
d(x;\complement_K^{\complement_K^P}) > r\rbrace.$ You hence just have to apply the previous
result to $\complement_K^P$.

\begin{flushright}
$\blacksquare$
\end{flushright}

So we have now an easier method for the calculation of the measure of a
Borel set: find a sequence of elements $M$ of known measures approaching by
above the Borel set (the best situation is when the limit notion $L$ can be
taken as a classical limit). Since $\mu_K$ is regular, the measure of the
Borel set is the limit of the measures of the elements of this sequence.
This method is even simpler in the case of closed sets.

\begin{proposition}
Let $K$ be a compact metric space. Then, for all closed sets $F$ of $K$,

\begin{equation*}
\mu_K(F) = \lim_{\delta \rightarrow 0} L_{\varepsilon \rightarrow 0} \frac{N(%
\overline{F}^{\delta};\varepsilon)}{N(K;\varepsilon)}. 
\end{equation*}
\end{proposition}

\vspace{2 mm}

\underline{\textbf{Demonstration:}} The density lemma enables to say that
the set $C$ of $\delta>0$ such as $\overline{F}^{\delta}$ is not in $M$ is
countable. So, for all $\delta \not\in C$, $\mu_K(\overline{F}^{\delta}) =
L_{\varepsilon \rightarrow 0} \frac{N(\overline{F}^{\delta};\varepsilon)}{%
N(K;\varepsilon)}$. Now, $\bigcap_{\delta > 0} \overline{F}^{\delta} = F$,
and this intersection is decreasing. Hence, one can choose a
decreasing sequence of $\delta$ converging to zero, with elements taken out
of $C$, say $(\delta_n)$, with $\bigcap_{n \in \mathbf{N}} \overline{F}%
^{\delta_n} = F$. Then:

\begin{equation*}
\mu_K(F) = \lim_{n\rightarrow +\infty} \mu_K(\overline{F}^{\delta_n}) =
\lim_{n\rightarrow +\infty}L_{\varepsilon \rightarrow 0} \frac{N(\overline{F}%
^{\delta_n};\varepsilon)}{N(K;\varepsilon)} = \lim_{\delta \rightarrow 0, \delta \neq 0}
L_{\varepsilon \rightarrow 0} \frac{N(\overline{F}^{\delta};\varepsilon)}{%
N(K;\varepsilon)}, 
\end{equation*}
since $\delta \mapsto L_{\varepsilon \rightarrow 0} \frac{N(\overline{F}%
^{\delta};\varepsilon)}{N(K;\varepsilon)}$ is decreasing.

\begin{flushright}
$\blacksquare$
\end{flushright}

This property of $\mu$ for closed enables us to state the following theorem,
using the regularity of our measure. This equivalent definition has a key
point: it does not use the algebra $M$.

\vspace{2 mm}

\begin{theorem}
Let $K$ be a compact metric space. For all Borel set $B$ of $K$, 
\begin{equation*}
\mu_K(B) = \sup_{ F \subseteq B, F closed} \lim_{\delta \rightarrow
0, \delta \neq 0} L_{\varepsilon \rightarrow 0} \frac{N(\overline{F}^{\delta};\varepsilon)}{%
N(K;\varepsilon)}. 
\end{equation*}
\end{theorem}

\vspace{2 mm}

The previous method is a bit heavy, despite its great generality. In
simplest cases, one should use the characteristic properties of $\mu_K$ ($%
\sigma$-additivity et invariance by partial isometries between opens) to
calculate measures of Borel sets. However, there is a difficult point:
with this technique, one shows that for all $\mu_K^L$, the measure of the
particular Borel set has a determined value, which does not depend on the
choice of the limit notion $L$. Hence, this method can not be used in all
situations. Indeed, $\mu_K^L$ sometimes depends on the choice of $L$.

\subsubsection{Compatibility with respect to product}

\vspace{2 mm}

\begin{proposition}
Let $A$ and $B$ be compacts. We endow $A\times B$ with the product
distance. Choose a $\varepsilon>0$. Then: $N(A\times B;\varepsilon) =
N(A;\varepsilon)N(B;\varepsilon)$.
\end{proposition}

\vspace{2 mm}

\underline{\textbf{Demonstration:}} If $R$ is a covering of $A$ by open
balls of radius $\varepsilon$, and if $R^{\prime }$ is such a covering of $B$%
, the family of products of a ball of $R$ and a ball of $R^{\prime }$ is a
covering of $A\times B$. Hence: $N(A\times B;\varepsilon) \leq
N(A;\varepsilon)N(B;\varepsilon)$.

Let then $\lbrace B_i\rbrace$ be a $\varepsilon$-minimal covering of $%
A\times B$ by open balls. Put for any subset $E$ of $A\times B$, $p_x(E) =
\lbrace y\in B / (x;y) \in E \rbrace$ where $x\in A$, and $q_y(E) = \lbrace
x\in A / (x;y) \in E \rbrace$ where $y\in B$. Then, for all $x$ in $A$, $P_x
= \lbrace p_x(B_i) / p_x(B_i) \neq \emptyset \rbrace$ is a covering of $B$
by $\varepsilon$-open balls, hence its cardinal is taller than $%
N(B;\varepsilon)$. Indeed, the image of a $\varepsilon$-ball for the product
distance under $p_x$ is a $\varepsilon$-ball if non-empty (a ball for the
product distance is a product of balls). In the same way, for all $y$ in $B$%
, $Q_y = \lbrace q_y(B_i) / q_y(B_i) \neq \emptyset \rbrace$ is a covering
of $A$ by $\varepsilon$-open balls, hence its cardinal is taller than $%
N(A;\varepsilon)$.

Then $\lbrace B_i\rbrace = \lbrace A^{\prime }\times B^{\prime }/ B^{\prime
}\in P_x, x \in A^{\prime }, A^{\prime }\in Q_y, y\in B \rbrace$ (show it
with a double inclusion). It is enough to show that the cardinal of that
last set is taller than $N(A;\varepsilon)N(B;\varepsilon)$.

Let us show that: choose a $y_0$ in $B$. Then $Q_{y_0}$ contains at least $%
N(A;\varepsilon)$ elements denoted $A^{\prime }_k$. Then, choose $%
N(A;\varepsilon)$ elements $(x_n)$ in different balls of $Q_{y_0}$ such that none of them is contained in two balls $A^{\prime }_k$. It is possible because if we could take $N
< N(A;\varepsilon)$ such elements, any more one would be
contained in an already chosen ball. Consequently, one could cover $A$ with
only $N$ balls, which is absurd. Then, for any $x_n$, we can choose at least 
$N(B;\varepsilon)$ different balls $B^{\prime }_l$ in $P_{x_n}$. Then,
letting $k$ and $l$ vary, we get a family of $N(A;\varepsilon)N(B;\varepsilon)$
different balls $A^{\prime }_k\times B^{\prime }_l$ in $\lbrace B_i\rbrace$.

\begin{flushright}
$\blacksquare$
\end{flushright}

\vspace{2 mm}

\begin{theorem}
Let $K$ and $K^{\prime }$ be non-empty compact spaces, and denote $\mu_Q$
the probability semi-canonically associated with any non-empty compact
metric space $Q$. Then $\mu_{K\times K^{\prime }}=\mu_K \otimes
\mu_{K^{\prime }}$, where $\mu_K \otimes \mu_{K^{\prime }}$ is the product
measure of $\mu_K$ and $\mu_{K^{\prime }}$, $K\times K^{\prime }$ being
endowed with the product distance.
\end{theorem}

\vspace{2 mm}

\underline{\textbf{Demonstration:}} Let $B$ be a Borel set of $K$, and $%
B^{\prime }$ a Borel set of $K^{\prime }$. Choose a $\varepsilon > 0$. The measures $\mu_K$ and $\mu_{K^{\prime }}$ being regular,
there are closed sets $F$ of $K$ and $F^{\prime }$ of $K^{\prime }$ such
that $F\subseteq B$, $F^{\prime }\subseteq B^{\prime }$, and $\mu_K(B) \geq
\mu_K(F)\geq \mu_K(B)-\varepsilon$, $\mu_K(B^{\prime }) \geq \mu_K(F^{\prime
})\geq \mu_K(B^{\prime })-\varepsilon$. Now, 
\begin{equation*}
\mu_{K\times K^{\prime }}(F\times F^{\prime }) = \lim_{\delta \rightarrow 0}
L_{\varepsilon \rightarrow 0} \frac{N(\overline{F\times F^{\prime }}%
^{\delta};\varepsilon)}{N(K\times K^{\prime };\varepsilon)} 
\end{equation*}
\begin{equation*}
= \lim_{\delta \rightarrow 0} L_{\varepsilon \rightarrow 0} \frac{N(%
\overline{F}^{\delta};\varepsilon)}{N(K;\varepsilon)}\frac{N(\overline{%
F^{\prime }}^{\delta};\varepsilon)}{N(K^{\prime };\varepsilon)} =
\mu_K(F)\mu_{K^{\prime }}(F^{\prime }). 
\end{equation*}
To see that just use the previous lemma and the fact that the Minkowski
sausage of a product is the product of the Minkowski sausages. Note that we use here the fact that $L$ is a morphism for the product law.

It follows that: $\mu_K(B)\mu_K(B^{\prime }) \geq \mu_{K\times K^{\prime
}}(F\times F^{\prime }) \geq \mu_K(B)\mu_K(B^{\prime})-2\varepsilon $. Then, taking $\varepsilon$ close to zero, we see that:
$$\mu_K(B)\mu_K(B^{\prime }) =
\sup_{F\subseteq B,F^{\prime }\subseteq B^{\prime }} \mu_{K\times K^{\prime
}}(F\times F^{\prime }) = \mu_{K\times K^{\prime }}(B\times B^{\prime }),$$
since $\mu_{K\times K^{\prime }}$ is regular.

The measures $\mu_{K\times K^{\prime }}$ and $\mu_{K}\otimes \mu_{K^{\prime
}}$ are then equal for products of Borel sets. Now the family of these
Borel sets generates the $\sigma$-algebra of Borel sets in the Dynkin systems sense. Thus these measures equal.

\begin{flushright}
$\blacksquare$
\end{flushright}

So we have shown a remarkable fact: the functor we are studying
commutes with finite products.

\subsubsection{Compatibility with conditioning}

If $F$ is a non-null closed of $K$, there are two probabilities on $F$: the
probability induced by that of $K$, $\frac{\mu_K}{\mu_K(F)}$, and the
measure associated to $F$, $\mu_F$. Actually, those probabilities do not
equal in whole generality.

Here is a counter example. Put $K = \lbrace \frac{1}{n}, n \in \mathbf{N}^*%
\rbrace \cup\lbrace 1 + \frac{1}{2n}, n \in \mathbf{N}^*\rbrace \cup \lbrace 0;1\rbrace$. We can easily show that $\mu_K = \frac{2}{3}\delta_0 + \frac{1}{3}\delta_1$. (By the way,
it illustrates the fact that $\mu_K$ gives some information about the
"density" of $K$.) Put $F = \lbrace 0;1\rbrace$ and $B = \lbrace 1\rbrace$.
One has: $\frac{\mu_K (B)}{\mu_K(F)} = \frac{1}{3}$, but $\mu_F (B) = \frac{%
1}{2}$.

Nevertheless, the probabilities $\mu_F$ and $\frac{\mu_K}{\mu_K(F)}$ are
linked in some way.

\vspace{2 mm}

\begin{proposition}
Let $K$ be a compact metric space, and $F$ a closed of $K$ with $\mu_K(F)
\neq 0$. Then, for any Borel set $B$ of $F$, $\lim_{\delta \rightarrow 0}
\mu_{\overline{F}^\delta}(B) = \frac{\mu_K(B)}{\mu_K(F)}$. Moreover, if $F
\in M_K$, for any Borel set $B$ of $F$, $\mu_F(B) = \frac{\mu_K(B)}{\mu_K(F)}
$.
\end{proposition}

\vspace{2 mm}

\underline{\textbf{Demonstration:}} Let us begin with the second part of
the proposition. Let $F \in M_K$ be a closed such that $\mu_K(F) \neq 0$. Then,
on one hand: $L_{\varepsilon \rightarrow 0}\frac{N(F;\varepsilon)}{%
N(K;\varepsilon)} = \mu_K(F) \neq 0$, and on the other hand, $\mu_K(Fr(F)) =
0$, hence $\overset{\circ}{F} \neq \emptyset$. Let then $B$ be a closed
positively included in $F$. One has:

\begin{equation*}
\mu_F(B) = \lim_{\delta \rightarrow 0}L_{\varepsilon \rightarrow 0} \frac{N(%
\overline{B}^{\delta}\cap F;\varepsilon)}{N(F;\varepsilon)} = \lim_{\delta
\rightarrow 0}L_{\varepsilon \rightarrow 0} \frac{N(\overline{B}%
^{\delta};\varepsilon)}{N(F;\varepsilon)} 
\end{equation*}
\begin{equation*}
= \lim_{\delta \rightarrow 0}L_{\varepsilon \rightarrow 0} \frac{N(\overline{%
B}^{\delta};\varepsilon)}{N(K;\varepsilon)} \frac{N(K;\varepsilon)}{%
N(F;\varepsilon)}= \frac{\mu_K(B)}{\mu_K(F)}. 
\end{equation*}

Indeed, in these equalities, for $\varepsilon$ small enough, the $N(.;\varepsilon)$ give the cardinalities
of minimal coverings in the sense of balls of $F$, but in sense of balls of $%
K$ too, since $B$ is positively included in $F$. Moreover, $L_{\varepsilon
\rightarrow 0}$ is an algebra morphism.

Now, denote $F_n = \overline{\lbrace x \in F / B(x; 1/n) \subseteq F\rbrace}$%
. The $F_n$ are closed positively included in $F$, and their union is
increasing and equalling $F$. Let then $B$ be any closed of $F$. $B =
\bigcup_n B\cap F_n$ is an increasing union of closed positively included in 
$F$. Thus: $\mu_F(B) = \lim_{n \rightarrow +\infty} \mu_F(B\cap F_n) =
\lim_{n \rightarrow +\infty} \frac{\mu_K(B\cap F_n)}{\mu_K(F)} = \frac{%
\mu_K(B)}{\mu_K(F)}.$ $\mu_F$ and $\frac{\mu_K}{\mu_K(F)}$ equal on closed sets.
Since they are measures on Borel sets, they are regular, and thus equal.

\vspace{2 mm}

Let us study now the case where $F$ is a closed not supposed to be in $M_K$.
Let $B$ be a Borel set of $F$. For all $\varepsilon > 0$, choose $\delta \in
(0;\varepsilon)$ such as $\overline{F}^{\delta} \in M_K$. $\frac{\mu_K(B)}{%
\mu_K(\overline{F}^{\delta})} = \mu_{\overline{F}^{\delta}}(B)$, following
from what precede. Take the limit when $\varepsilon$ approaches zero: $%
\lim_{\delta \rightarrow 0} \mu_K(\overline{F}^{\delta}) = \mu_K(F)$ since 
$F$ is closed. Therefore: $\lim_{\delta \rightarrow 0} \mu_{\overline{F}%
^\delta}(B) = \frac{\mu_K(B)}{\mu_K(F)}$.

\begin{flushright}
$\blacksquare$
\end{flushright}

\subsection{Definitive version of the theorem. Notion of metrically
compatible probability of a compact metric space}

Now, we give stronger properties about the invariance properties of the probabilities we
built. We need the following definition to reach this aim.

\vspace{2 mm}

\begin{definition}
Let $E$ and $F$ be metric spaces. We say that a map $%
f : E \rightarrow F$ such that for all $x$ and $y$ of $E$, $d(f(x);f(y))
\geq d(x;y)$ is $1$-coercive, and that a map $g : E \rightarrow F$ such that for all $x$ and $y$ of $E$, $d(f(x);f(y))\leq d(x;y)$ is $1$-lipschitz.
\end{definition}

\vspace{2 mm}

\begin{proposition}
$1$-lipschitz maps between metric spaces are continuous, and $1$-coercive maps are injective. If $E$ is compact, any $1$-coercive map $f : E \rightarrow F$ is continuous if and only if $f(E)$ is a compact subset of $F$, and $f^{-1}: f(E)\rightarrow E$ is $1$-lipschitz. $f$ is consequently a homeomorphism onto its image under these hypothesis. Moreover, if $g : E \rightarrow F$ is $1$-lipschitz and injective, $g^{-1}:g(E)\rightarrow E$ is $1$-coercive with $g(E)$ compact.
\end{proposition}

\vspace{2 mm}

\underline{\textbf{Demonstration:}} The first two affirmations and the last one are obvious. Let then $f : E \rightarrow F$ be a $1$-coercive map. If $f$ is continuous, $f(E)$ is compact. Assume now that $f(E)$ is compact. Then, $f^{-1}: f(E)\rightarrow E$ is obviously $1$-lipschitz, and for any closed $C \subset f(E)$, $C$ is compact, and $f^{-1}(C)$ is compact, hence a closed of $E$. Thus, $f$ is continuous.\begin{flushright}
$\blacksquare$
\end{flushright}

We use properties shown in the previous part in order to prove inequalities
implying the invariance by partial isometries between opens.

\begin{proposition}
Let $K$ be a compact metric space and $\mu$ our probability on $K$. Let $U$
and $V$ be opens of $K$ such that there is a $1$-coercive continuous map $f$ from $U$
to $V$, not assumed to be surjective. Then: $\mu (U) \leq \mu (V)$. Moreover, if $U$ and $V$ are opens of $K$ such that there is a $1$-lipschitz map $g$ from $U$ onto $V$, not necessarily injective, then $\mu (U) \geq \mu (V)$.
\end{proposition}

\vspace{2 mm}

\underline{\textbf{Demonstration:}} Let $F$ be a closed set included in $U$, and $G = f(F)$. These inclusions are positive. Since $f$ is $1$-coercive, $G$ is a compact
set positively included in $V$, and if $(B_i)$ is a $\varepsilon$-minimal
covering of $G$, $(f^{-1}(B_i))_i$ is a covering of $F$. Moreover, since the
centres $c_i$ of the balls $B_i$ are in the image of $F$ by $f$, for all $i$%
, $f^{-1}(c_i)$ is non-empty, and is a singleton. Then, denote $f^{-1}(c_i) = \lbrace b_i\rbrace$. Since $f$ is $1$%
-coercive, for all $i$, $B(b_i; \varepsilon) \supseteq f^{-1}(B_i)$. Thus $%
N(F;\varepsilon) \leq N(G;\varepsilon)$ for all $\varepsilon > 0$.

Now, if $F$ belongs to $M_K$, $\mu (F) = L_{\varepsilon \rightarrow 0}\frac{N(F;\varepsilon)}{N(K;\varepsilon)} \leq L_{\varepsilon \rightarrow 0}\frac{N(G;\varepsilon)}{N(K;\varepsilon)}\leq \mu (G) \leq \mu (V)$ since $G$ is closed.

And with the density lemma, any open set is an increasing denumerable union of sets of $M_K$, hence $\mu (U) \leq \mu (V)$.

Let now $g : U \rightarrow V$ be a $1$-lipschitz map between opens of $K$. Let $G$ be any compact set of $V$. $f^{-1}(G) = F$ is a compact of $U$. If $(B_i)$ is a $\varepsilon$-minimal covering of $G$, $(f^{-1}(B_i))_i$ is a covering of $F$. Denoting $c_i$ the antecedent of the centre of $B_i$ by $f$ for all $i$, we have $B(b_i; \varepsilon) \subset f^{-1}(B_i)$ for all $i$. Hence $%
N(F;\varepsilon) \geq N(G;\varepsilon)$ for all $\varepsilon > 0$.

Then, if $G$ is in $M_K$, $\mu (G) = L_{\varepsilon \rightarrow 0}\frac{N(G;\varepsilon)}{N(K;\varepsilon)} \leq L_{\varepsilon \rightarrow 0}\frac{N(F;\varepsilon)}{N(K;\varepsilon)} \leq \mu (F) \leq \mu (U)$.

Now, with the density lemma, any open set is an increasing denumerable union of sets of $M_K$, hence $\mu (V) \leq \mu (U)$.\begin{flushright}
$\blacksquare$
\end{flushright}

Consequently, we put the following definition:

\begin{definition}
Let $E$ be a metric space. A metrically compatible measure of $E$ is a
measure $m$ defined on Borel sets such that:

- for any couple of opens $(U;V)$
of $E$ such that there is a continuous $1$-coercive map $\varphi : U
\rightarrow V$, $m(U)\leq m(V)$,

- for any couple of opens $(U;V)$
of $E$ such that there is a $1$-lipschitz map $\varphi : U
\rightarrow V$, $m(U)\geq m(V)$.
\end{definition}

\vspace{2 mm}

\begin{proposition}
If $m$ is a metrically compatible measure in a metric space $E$, for any couple of
partially isometric opens $(U;V)$, $m(U) = m(V)$.
\end{proposition}

\vspace{2 mm}

\begin{definition}
Let $E$ be a metric space. An isometrically compatible measure of $E$ is a
measure $m$ defined on Borel sets such that for any couple of opens $(U;V)$
of $E$ such that there is a partial isometry $\varphi : U
\rightarrow V$, $m(U)= m(V)$.
\end{definition}

Obviously, metrically compatible measures are isometrically compatible.

\vspace{2 mm}

For the following proposition, beware to the distinction between isometries and isometric maps. The latter are not supposed to be surjective.

\vspace{2 mm}

\begin{proposition}
Let $m$ be a metrically compatible $\sigma$-finite measure on a metric space $E$.

- If $\varphi : U \rightarrow V$ is continuous and $1$-coercive between opens, and if $P$ is any Borel set of $U$, $m(P)\leq m(\varphi (P))$,

- If $\varphi : U \rightarrow V$ is $1$-lipschitz between opens, and if $P$ is any Borel set of $U$, $m(P)\geq m(\varphi (P))$,

- If $\varphi : U \rightarrow V$ is an isometric map between opens, and if $P$ is any Borel set of $U$, $m(P)= m(\varphi (P))$.

Moreover, if $m$ is just an isometrically compatible $\sigma$-finite measure on $E$, and if $\varphi : U \rightarrow V$ is an isometry between opens, any Borel subset $P$ of $U$ is such that: $m(P)= m(\varphi (P))$.
\end{proposition}

\vspace{2 mm}

\underline{\textbf{Demonstration:}} $m$ is a regular measure in both situations. Thus, for any partial continuous map $\varphi : U \rightarrow V$, all Borel set $P \subset U$ and all $\varepsilon > 0$, there are opens $O$ and $O'$ such that $P \subset O$, $f(P) \subset O'$, and $m(O) < m(P) + \varepsilon$, $m(O') < m(f(P)) + \varepsilon$. Then $\Omega = U\cap O \cap \varphi^{-1}(O')$ is open, contains $P$, is contained in $O$, and is such that $\varphi (\Omega) \subset O'$. Moreover, $m(\Omega) < m(P) + \varepsilon$ and $m(\varphi (\Omega)) < m(\varphi (P)) + \varepsilon$.

Now, if $m$ is metrically compatible and $\varphi$ $1$-coercive, $m(\Omega)\leq m(\varphi (\Omega))$, which leads to the conclusion.

The other cases can be treated similarly.\begin{flushright}
$\blacksquare$
\end{flushright}

The set of metrically compatible measures and this of isometrically compatible ones of a non-empty compact metric
space are non-empty bounded convex sets of the space of Radon measures. (Both are non-empty thanks to our main theorem.) There are many interesting facts
following this remark: indeed, these spaces are compact in weak-$*$ topology, which enables using
interesting theorems, such as the Krein-Milman one, or powerful techniques linked with
compactness. Maybe one could use techniques or theorems from ergodic theory
to show some more results.

Note that these spaces do not necessarily equal: for example, look at $K = \lbrace \frac{1}{n}, n \in \mathbf{N}^*\rbrace \cup \lbrace 2+\frac{1}{2n}, n \in \mathbf{N}^*\rbrace \cup\lbrace 0;2\rbrace$. The Dirac measure in $2$ is isometrically compatible but not metrically compatible.

\begin{proposition}
Let $K$ and $K^{\prime }$ be compact metric spaces, and $f$ be an isometric
map between $K$ and $K^{\prime }$. Then, for all
Borel set $P$ of $K^{\prime }$, $$\mu_{K}(f^{-1}(P)) \geq \mu_{K^{\prime
}}(P) .$$
\end{proposition}

\vspace{2 mm}

\underline{\textbf{Demonstration:}} Let $Q$ be the compact image of $K$ by $%
f$. One has: 
\begin{equation*}
\mu_{K}(f^{-1}(U)) = \mu_{Q}(U) = L_{\varepsilon \rightarrow 0} \frac{%
N(U;\varepsilon)}{N(Q;\varepsilon)} 
\end{equation*}
\begin{equation*}
\geq L_{\varepsilon \rightarrow 0} \frac{N(U;\varepsilon)}{N(K^{\prime
};\varepsilon)} \geq \mu_{K^{\prime }}(U) 
\end{equation*}
for all opens $U \in M_{Q}$. It follows using the density lemma that for all
opens $O$ of $K^{\prime }$, $\mu_K(f^{-1}(O)) \geq \mu_{K^{\prime }}(O)$,
and then, using the regularity of $\mu_{K}$ and $\mu_{K^{\prime }}$, that for all
Borel sets $P$ of $K^{\prime }$, $\mu_{K}(f^{-1}(P)) \geq \mu_{K^{\prime
}}(P)$.

\begin{flushright}
$\blacksquare$
\end{flushright}

The following theorem sums up our results.

\vspace{2 mm}

\begin{theorem}[Geometry and probability]
Let $L$ be a limit notion in $0$ of the space of bounded functions on $%
(0;+\infty )$. On all compact metric spaces $K$, there is a semi-canonical
metrically compatible probability $\mu_K$ depending only on the choice of 
$L$, such that:

- if $K$ and $K^{\prime }$ are isometric compacts, $\mu_{K}$ and $%
\mu_{K^{\prime }}$ are image from each other via this isometry,

- $\mu_{K\times K^{\prime }} = \mu_{K} \otimes \mu_{K^{\prime }}$, where $%
K\times K^{\prime }$ is endowed with the product distance,

- if $F$ is a closed of $K$ such that $\mu_K(F) \neq 0$, for all Borel sets $%
B$ of $F$, $\lim_{\delta \rightarrow 0} \mu_{\overline{F}^\delta}(B) = \frac{%
\mu_K(B)}{\mu_K(F)}$,

- if $f$ is an isometric map between $K$ and $K^{\prime }$, $f_* \mu_{K}
\geq \mu_{K^{\prime }}$.

We call those probabilities the metrically compatible probabilities
canonically associated with the limit notion $L$.
\end{theorem}

\vspace{2 mm}

So, choosing a limit notion enables us to endow compact metric spaces with a
probability measure satisfying the previous properties. Formalising a bit
further, we can see that the choice of a limit notion $L$ determines a
functor from the category of compact metric spaces and isometric maps to the
category of compact metric spaces endowed with a metrically compatible
probability and isometric maps such that if $f: (K;\mu_K) \rightarrow (K';\mu_{K'})$, $f_{*}\mu_K \geq \mu_{K'}$. Moreover, this
functor commutes with finite products, and its composition with the
probability forgetful functor is the identity one.

\vspace{2 mm}

The fact that the choice of a limit notion determines a canonical way to naturally endow each compact metric space with a metrically compatible probability is quite encouraging for the study of these spaces. Indeed, some interesting constructions generalizing these usually made in manifolds can be done choosing a Radon measure. See for example the case of the Villani-Ricci curvature (\cite{vill1} and \cite{vill2}). With our theorem, such constructions do not need any more such a choice, since the measure comes naturally from the topology of the space (once a limit notion is given). One can still endow the compact space with another measure (even not a metrically compatible one), but in case of uniqueness of the metrically compatible probability, this choice is the most natural one in the perspective of a generalization of the Geometry of manifolds.

\section{The question of uniqueness}

In order to give some properties about uniqueness of metrically
compatible measures, we set the following definitions.

\subsection{Homogeneity of a metric space}

\begin{definition}
Let $M$ be a metric space. We define the following relation on points of $M$: $x \leftrightarrow y$ $\Leftrightarrow$ there is a $\varepsilon>0$ and an
isometry from $B(x;\varepsilon)$ onto $B(y;\varepsilon)$. $\leftrightarrow $
is an equivalence relation, called the homogeneity equivalence of $K$.
\end{definition}

\vspace{2 mm}

\begin{definition}
With the previous notations, $M$ is said homogeneous if $\leftrightarrow$
has only one equivalence class. Otherwise, the equivalence classes of $\leftrightarrow$ are called the homogeneous components of $M$.
\end{definition}

\vspace{2 mm}

Now, we define an other notion linked to the previous one.

\begin{definition}
Let $M$ be a metric space. We define the following relation on points of $M$: $x \hookrightarrow y$ $\Leftrightarrow$ there is a $\varepsilon>0$ and an
isometric map from $B(x;\varepsilon)$ to $B(y;\varepsilon)$. $\hookrightarrow $
is a pre-order relation call the homogeneity pre-order.
\end{definition}

\vspace{2 mm}

\begin{proposition}
Let $M$ be a compact metric space. Then the kernel of the relation $\hookrightarrow$ is the relation $\leftrightarrow$.
\end{proposition}

\vspace{2 mm}

\underline{\textbf{Demonstration:}} Let $x$ and $y$ be elements of $M$ with $x \hookrightarrow y$ and $y \hookrightarrow x$. Then, one can find a $\varepsilon > 0$ and isometric maps $f : B(x;\varepsilon) \rightarrow B(y;\varepsilon) $ and $g : B(y;\varepsilon) \rightarrow B(x;\varepsilon) $, and $g\circ f$ is an isometric map from $B(x;\varepsilon)$ to itself. Hence, its restriction $h$ to $ K = \overline{B(x;\varepsilon /2)}$ is isometric too, and goes into $K$. So $h(K)\subset K$ is compact, and if $h$ is not surjective, there is a $x \in K$ such that $x \not\in h(K)$. Thus $d(x; h(K)) = \alpha >0$. Consider now the sequence $(h^n(x))_n$ of iterates of $h$ applied on $x$. For any integers $n$ and $p$ with $p>n$, $d(h^n(x);h^p(x)) = d(x;h^{p-n}(x))> \alpha$. Hence $(h^n(x))_n$ cannot have any limit point. But $K$ is compact : it is absurd.

Thus $h$ is surjective. And since it is an isometry, the restriction of $h$ to $B(x;\varepsilon /2)$ is surjective onto $B(x;\varepsilon /2)$, and $f$ is an isometry between $B(x;\varepsilon /2)$ and $B(y;\varepsilon /2)$. Hence $x\leftrightarrow y$.
\begin{flushright}
$\blacksquare$
\end{flushright}

The following lemma will be useful later.

\begin{lemma}
Let $E$ be a non-empty metric space. Let $x\in E$, and $[x]$ its class for the homogeneity equivalence. Put: 
\begin{equation*}
H_n (x) = \lbrace y \in E / \mbox{ there is an isometry between } B(x;1/n) %
\mbox{ and } B(y;1/n)\rbrace. 
\end{equation*}
Then, for all $y$ and $z$ in $H_n (x)$, $B(y;1/n)$ and $B(z;1/n)$ are isometric, and $[x] = \bigcup_n H_n (x)$. Moreover, if $E$ is compact, $H_n (x)$ is closed for all $n$, which implies that $[x]$ is a Borel set.
\end{lemma}

\vspace{2 mm}

\underline{\textbf{Demonstration:}} The first and second affirmations are obvious. Now, let $(y_k)$ be sequence of $E$ converging to $y \in E$, and $\Phi_k : B(x;1/n) \rightarrow B(y_k;1/n)$ be an isometry. We will use an argument similar to that used to prove the Ascoli's Theorem. $E$ is separable, hence we can choose a denumerable family $(p_i)$ dense in $B(x;1/n)$. Then, for $i=0$, $\Phi_k (p_0)$ has a limit point in $\overline{B(y;1/n)}$, that is, it has a converging subsequence $\Phi_{\varphi_0(k)} (p_0)$. Then, for $i=1$, $\Phi_{\varphi_0(k)} (p_1)$ admits a converging subsequence $\Phi_{\varphi_0(\varphi_1(k))} (p_1)$.

Continuing so on, we recursively build strictly increasing sequences of integers $\varphi_j$. Now, with the Cantor diagonal process, taking : 
$$\phi (k) = (\bigcirc_{j=0}^{k}\varphi_j)(k),$$
we define a strictly increasing sequence of integers such that $\Phi_{\phi(k)} (p_i)$ converges for all $i$. Using these limits, we define a map $f' : \lbrace p_i, i \in \mathbf(N) \rbrace \rightarrow \overline{B(y;1/n)}$. $\Phi_k$ being an isometry for all $k$, $f'$ is an isometric map, and it is uniformly continuous. Hence, it can be continuously prolonged to $\overline{B(x;1/n)}$ into an map $f : \overline{B(x;1/n)} \rightarrow \overline{B(y;1/n)}$. One can easily show that $f$ is an isometric map, and that its restriction $g$ to $B(x;1/n)$ is an isometric map into $B(y;1/n)$.

Now, let $z$ be any element of $B(y;1/n)$. $d(y;z) < 1/n$, hence for $k$ tall enough, say for $k \geq K$, $z$ belongs to $B(y_{\phi(k)};1/n-\frac{1/n-d(y;z)}{2})$. Hence $\Phi_{\phi(k)}^{-1} (z)$ is a defined element of $\overline{B(x;1/n-\frac{1/n-d(y;z)}{2})}$, a compact set. Hence, this sequence has a limit point $u$. But, since for all $k \geq K$, $\Phi_{\phi(k)}(\Phi_{\phi(k)}^{-1}(z))=z$, 

$$d(f(u);z) \leq $$
$$d(f(u);\Phi_{\phi(k)}(u)) + d(\Phi_{\phi(k)}(u); \Phi_{\phi(k)}(\Phi_{\phi(k)}^{-1}(z)))+d(\Phi_{\phi(k)}(\Phi_{\phi(k)}^{-1}(z));z) =$$
$$d(f(u);\Phi_{\phi(k)}(u)) + d(u; \Phi_{\phi(k)}^{-1}(z)).$$

Thus, to conclude that $f(u) = z$, and that $g$ is an isometry, which shows that $y \in H_n (x)$, it is sufficient to prove that $\Phi_{\phi(k)}$ simply converges to $f$. 

Let $u \in \overline{B(x;1/n)}$, and choose any $\varepsilon > 0$. There is a $p_i$ such that $d(u;p_i) < \varepsilon$, and a $K$ such that $d(\Phi_{\phi(k)} (p_i),f(p_i))< \varepsilon$ for all $k>K$. Thus : 
$$d(\Phi_{\phi(k)} (u),f(u))\leq $$
$$d(\Phi_{\phi(k)} (u);\Phi_{\phi(k)} (p_i))+d(\Phi_{\phi(k)} (p_i);f(p_i))+d(f(p_i);f(u)) < $$
$$\varepsilon + 2d(p_i;u)=3\varepsilon.$$

This shows that $\Phi_{\phi(k)}$ simply converges to $f$, which ends the proof.\begin{flushright}
$\blacksquare$
\end{flushright}

\subsection{A sufficient condition for uniqueness of metrically
compatible probabilities}

We use here an adapted version of the Christensen theorem shown in \cite%
{Chris}.

\begin{definition}
Let $m$ be a measure on a metric space $M.$ $m$ is said uniform if there is
an increasing family $(M_n)$ of Borel sets of $M$ whose union is $M$ and a
sequence $(r_n)$ of positive real numbers, such that for all natural numbers 
$n$, and all $x$ and $y$ in $M_n$, and for all positive $r^{\prime }<r_n$, $%
m(B(x;r^{\prime })) = m(B(y;r^{\prime }))$.
\end{definition}

\vspace{2 mm}

Note that in this definition, there is no need for the balls $B(x;r^{\prime
})$ to be included in $M_n$ for $r^{\prime }< r_n$.

\vspace{2 mm}

With the definition above and the last lemma of the previous paragraph, we easily see that isometrically compatible probabilities are uniform measures when restricted to the homogeneous components of any compact metric space.

\begin{theorem}[Christensen revisited]
Let $P\neq 0$ and $P^{\prime }$ be uniform non-negative finite measures on $M$
defined on Borel sets, where $M$ is a pre-compact metric space. Then, there
is a $l \in \mathbf{R}$ such that $P^{\prime }= lP$. Consequently, if $P$
and $P^{\prime }$ are probabilities, they equal.
\end{theorem}

\vspace{2 mm}

\underline{\textbf{Demonstration:}} First of all, note that if there is a $%
x \in M$ such that $P^{\prime }(\lbrace x\rbrace)\neq 0$, $P^{\prime }$ is a
multiple of the counting measure on a finite space (use finiteness of $%
P^{\prime }$). $P$ being uniform too, $P^{\prime }$ is a multiple of $P$.

Assume now that $P^{\prime }(\lbrace x\rbrace) = 0$ for all $x$ in $M$,
which implies that: 
\begin{equation*}
\lim_{r \rightarrow 0} P^{\prime }(B(x;r)) = 0. 
\end{equation*}

\vspace{2 mm}

$P$ and $P^{\prime }$ are uniform, hence, there are increasing sequences $%
(M_n)$ and $(M^{\prime }_n)$ with union $M$ and sequences $(r_n)$ and $%
(r^{\prime }_n)$ satisfying the conditions of the definition. Then $%
M^{\prime \prime }_n = M_n \cap M^{\prime }_n$ and $r^{\prime \prime }_n =
\min (r_n;r^{\prime }_n)$ satisfy those conditions too for the uniform
measures $P$ and $P^{\prime }$. Now, choose a natural number $n$, and define the function $(x;y)\mapsto K_r^n(x;y)$ on ${M}^2$ for $r\in
(0;r^{\prime \prime }_n)$ by:

$K_r^n(x;y) = \frac{1}{c_r(P)}$ if $d(x;y) < r$ and $x \in M^{\prime \prime
}_n$, and

$K_r^n(x;y) = 0$ otherwise,

where $c_r(P) = P(B(x;r))$ does not depend on the choice of $x$. It is
non-zero for $n$ large enough. Indeed, if $c_r(P) = 0$, one can cover $%
M^{\prime \prime }_n$ by a finite number of null balls $B(z;r)$
(pre-compactness), and $M^{\prime \prime }_n$ is a null set. If it is true
for any $n$, $M = \bigcup_n M^{\prime \prime }_n$ is a null set, which is
absurd since $P\neq 0$. By the way, this shows that for $n$ tall enough, $M^{\prime \prime }_n$ is not a $P$-null set.

Then, 
\begin{equation*}
\int_M K_r^n(x;y)dP(y) = \mathbf{1}_{M^{\prime \prime }_n}(x), 
\end{equation*}

and:

\begin{equation*}
\int_M K_r^n(x;y)dP(x) = \frac{P(B(y;r)\cap M^{\prime \prime }_n)}{c_r(P)}. 
\end{equation*}

Now, let $f$ be a continuous function with compact support in $M^{\prime \prime}_n$, and define $(K_r^n f)(x) = \int_{M}
K_r^n(x;y)f(y)dP(y).$ One has:

\begin{equation*}
\vert (K_r^nf)(x) - f(x)\vert \leq \int_{M^{\prime \prime }_n}
K_r^n(x;y)\vert f(y) - f(x) \vert dP(y) \leq W_r(f), 
\end{equation*}

where 
\begin{equation*}
W_r(f) = \sup \lbrace\vert f(y) - f(x) \vert, d(x;y)<r\rbrace. 
\end{equation*}

Since $f$ is continuous and has compact support, it is uniformly continuous,
and $\lim_{r \rightarrow 0} W_r(f) = 0.$ Hence $\lim_{r \rightarrow 0}
(K_r^nf)(x) = f(x)$ for all $x$ in $M$. But one shows easily that $K_r^nf$
has compact support for $r$ small enough, and is uniformly bounded in $r$. Using the Lebesgue
Theorem, we get: 
\begin{equation*}
\lim_{r \rightarrow 0} \int_{M} (K_r^nf)(x)dP^{\prime }(x) = \int_{M}
f(x)dP^{\prime }(x) 
\end{equation*}
for all $n$.

Put then $c_r(P^{\prime }) =P^{\prime }(B(x;r))$, a constant in $x$. With
Fubini Theorem:

\begin{equation*}
\int_{M} (K_r^nf)(x)dP^{\prime }(x) = \int_{M} \int_{M}
K_r^n(x;y)f(y)dP(y)dP^{\prime }(x) = 
\end{equation*}
\begin{equation*}
\int_{M} \int_{M^{\prime \prime}_n} \frac{1}{c_r(P)}\textbf{1}_{B(x;r)}(y)f(y)dP(y)dP^{\prime }(x) = 
\end{equation*}
\begin{equation*}
\frac{1}{c_r(P)}\int_{M^{\prime \prime}_n} (\int_{M}\textbf{1}_{B(y;r)}(x)dP^{\prime }(x))f(y)dP(y) 
\end{equation*}
\begin{equation*}
= \frac{c_r(P^{\prime })}{c_r(P)} \int_{M^{\prime \prime }_n} f(y)dP(y) = \frac{c_r(P^{\prime })}{c_r(P)} \int_{M} f(y)dP(y). 
\end{equation*}

Choosing $f$ such that $\int_{M} f(y)dP(y) \neq 0$ (it is
possible  for $n$ tall enough when $M^{\prime \prime}_n$ is not $P$-null, since $P$ is regular), it
shows that $\lim_{r \rightarrow 0} \frac{c_r(P^{\prime })}{c_r(P)}$ exists.
Moreover, that limit does not depend on $n$. Denote it $l$. Take now the limit
as $r$ approaches $0$:

\begin{equation*}
\int_{M^{\prime \prime }_n} f(x)dP^{\prime }(x) = l \int_{M^{\prime \prime
}_n} f(y)dP(y). 
\end{equation*}

Since $l$ is fixed independently from $f$, it follows that $lP$ and $%
P^{\prime }$ equal on all $M^{\prime \prime }_n$, which implies that they
equal on $M$.

\begin{flushright}
$\blacksquare$
\end{flushright}

\vspace{2 mm}

Once this result is known, we get the following theorem.

\begin{theorem}
Let $K$ be a pre-compact homogeneous metric space, and $P \neq 0$ and $P'$ be isometrically compatible measures on $K$. Then, there is a $l \in \mathbf{R}$ such that $P'=lP$. Consequently, if $P$ and $P'$ are probabilities, they equal.
\end{theorem}

\vspace{2 mm}

\underline{\textbf{Demonstration:}} With the Christensen theorem, it is enough to prove that $P$ and $P'$ are uniform. Just apply the lemma shown in the paragraph about homogeneity: $H$ is the denumerable union of $(H_n)$. For all $n$ and all $x$ and $y$
in $H_n$, there is an isometry from $B(x;1/n)$ to $B(y;1/n)$, hence for all
non-negative $r^{\prime }<1/n$, $P(B(x;r^{\prime })) = P(B(y;r^{\prime }))$,
and the same holds for $P^{\prime }$. 
\begin{flushright}
$\blacksquare$
\end{flushright}

As a consequence, our theorems show that on any compact homogeneous space, there is a unique isometrically compatible probability. Moreover, any homogeneous component of a compact set is a pre-compact set. Thus, if it admits an isometrically compatible probability, this one is unique.

An interesting theorem describing the space of isometrically compatible probabilities under a suitable assumption follows.

\vspace{2 mm}

\begin{theorem}
Let $K$ be a compact metric space, and $(K_i)_i$ its homogeneous components. Assume that homogeneous components are in countable cardinal. For all $i$, if $K_i$ admits a non-zero non-negative finite isometrically compatible Borel measure, denote $\mu_i$ its unique isometrically compatible probability, and denote $I$ the set of such $i$. Then, $\mu$ is an isometrically compatible probability of $K$ if and only if there is a family of positive real numbers $(\lambda_i)_{i\in I}$ such that $\sum_i \lambda_i = 1$ and $\mu = \sum_i \lambda_i \mu_i$. Consequently, the space of isometrically compatible probabilities of $K$ is a compact convex set (for $*$-weak topology) whose extreme points are isometrically compatible probabilities supported on homogeneous components of $K$ admitting an isometrically compatible probability.
\end{theorem}

\vspace{2 mm}

\underline{\textbf{Demonstration:}} The formula $\mu = \sum_i \lambda_i \mu_i$ defines an isometrically compatible probability of $K$. Indeed, let $f : O \rightarrow O'$ be any isometry between opens of $K$. $O = \coprod_i O\cap K_i$, and for each $i$, $O\cap K_i$ is a Borel set. Hence $\mu_i (O\cap K_i) = \mu_i (f(O\cap K_i))$, and $O' = \coprod_i f(O\cap K_i)$. Then, denoting $J = \lbrace i \in I / \lambda_i \neq 0\rbrace$, $J$ is countable and $\mu (O) = \sum_{i \in J} \lambda_i \mu_i(O\cap K_i) = \sum_{i \in J} \lambda_i \mu_i(f(O\cap K_i))$. On the other hand, $\mu (O') = \sum_{i \in J} \lambda_i \mu_i (\coprod_i f(O\cap K_i))$. But for any $x \in O\cap K_i$, $f$ can be restricted to define an isometry between a ball centred in $x$ onto another centred in $f(x)$. Hence $f(x) \in K_i$, and $f(O\cap K_i) \subset K_i$. So : $\mu (O') = \sum_i \lambda_{i \in J} \mu_i (f(O\cap K_i)) = \mu (O)$.

Let now $\mu$ be any isometrically compatible probability of $K$. Denote $\mu\vert_{K_i}$ its restriction to each homogeneous component $K_i$. Looking at the definition of homogeneous components, we see that $\mu\vert_{K_i}$ is isometrically compatible. Hence, $\mu\vert_{K_i} = l \mu_i$ for some $l$, and so, $\mu\vert_{K_i} = \mu(K_i) \mu_i$. Denote $\mu ' (B) = \sum_i \mu(K_i) \mu_i (B\cap K_i)$ for all Borel sets $B$. $\mu '$ is obviously a Borel measure, and for all Borel set $B$, $\mu ' (B) = \sum_i \mu\vert_{K_i}(B\cap K_i) =\sum_i \mu (B\cap K_i) = \mu (B)$, so $\mu = \mu '$. Finally, note that $\mu (K) = 1 = \sum_i \mu(K_i)$. One just have to take $\lambda_i = \mu(K_i)$.
\begin{flushright}
$\blacksquare$
\end{flushright}

The uniqueness result can be seriously strengthened for metrically compatibles probabilities using the homogeneity pre-order. We even show this complementary result for more general measures than just the metrically compatible ones.

\vspace{2 mm}

\begin{definition}
Let $E$ be a metric space, and $m$ be a Borel measure on $E$. If for any isometric map $f : O \rightarrow O'$ defined between open sets of $E$ ($f$ does not need to be surjective), $m(f(O)) = m(O)$, $m$ is said strongly isometrically compatible.
\end{definition}

Note that metrically compatible probabilities are strongly isometrically compatible, and that the latter are isometrically compatible. Finding spaces and measures showing that the reverses are false is an interesting exercise.
 
\vspace{2 mm}

\begin{theorem}
Let $K$ be a compact metric space, such that maximal elements of the homogeneity pre-order relation of $K$ are all homogeneously equivalent. Denote $K_m$ the greatest homogeneous component for the homogeneity pre-order, and assume that it has a non-empty interior. Then, there is only one strongly isometrically compatible probability on $K$.
\end{theorem}

\vspace{2 mm}

\underline{\textbf{Demonstration:}} Let $\mu$ and $\mu '$ be strongly isometrically compatible probabilities on $K$. They are uniform measures on $K_m$, hence, there is a $l > 0$ such that $\mu = l \mu '$ (assuming that $\mu '(K_m) \neq 0$).

Now, let $x$ be in any homogeneous component $K_i$ of $K$, and $y$ and $\alpha > 0$ such that $B(y;\alpha)$ is contained in $K_m$. There is a $\varepsilon_x \in (0;\alpha )$ such that $B(x;\varepsilon_x)$ can be mapped isometrically into $B(y;\varepsilon_x)$ by, say, $f_x$.

And the family of $B(x;\varepsilon_x)$ covers $K$ which is compact. Hence, $K$ can be covered by a finite family of such balls, say $(B(x_i;\varepsilon_{x_i}))$.

Let now $B$ be any Borel set of $K$. One has : $\mu (B) = \sum_i \mu (B\cap B(x_i;\varepsilon_{x_i})) = \sum_i \mu (f_{x_i}(B\cap B(x_i;\varepsilon_{x_i}))) = \sum_i l\mu ' (f_{x_i}(B\cap B(x_i;\varepsilon_{x_i}))) = l\sum_i \mu '(B\cap B(x_i;\varepsilon_{x_i})) = l \mu '(B)$. Thus $\mu = l \mu '$, and so $\mu = \mu'$ since they are both probabilities.

Now, with the same kind of proof, one shows that $\mu ' = 0$ if $\mu '(K_m) = 0$, justifying our first assumption.\begin{flushright}
$\blacksquare$
\end{flushright}

\subsection{Notion of intrinsically null sets}

We develop in this section a notion of null sets for all strongly isometrically compatible measures. This notion can reveal itself useful in some situations, even if our uniqueness theorems can apply.

\begin{definition}
Let $E$ be a compact metric space. A Borel set $P$ of $E$ is said strongly isometrically
intrinsically null if for any strongly isometrically compatible probability $m$ on $E$, $m(P) = 0$.
\end{definition}

\vspace{2 mm}

Definitions of similar notions using isometrically or metrically compatible probabilities exist, but we will not need them in this paper.

\vspace{2 mm}

Now, we give a sufficient condition for a Borel set to be strongly isometrically intrinsically null.

\begin{proposition}
Let $E$ be a metric space. Denote $\upsilon (E)$ the family of Borel sets $P$ of $E$ such that for any natural number $N$, there are an open $O$
containing $P$ and a finite family $(\phi_i)_{i \leq N}$ of isometric maps
from $O$ to $E$ such that the $(\phi_i (P))_{i \leq N}$ are mutually
disjoint.

Denote now $\Upsilon (E)$ the smallest family of subsets of $E$ containing $\upsilon (E)$ and invariant under countable unions. The
elements of $\Upsilon (E)$ are strongly isometrically intrinsically null.
\end{proposition}

\vspace{2 mm}

\underline{\textbf{Demonstration:}} It is enough to show the result for sets of $\upsilon (E)$. Indeed, elements of $\Upsilon (E)$ are just countable unions of elements of $\upsilon (E)$. Let $P$ be such a set. Then, for any
natural number $N$, there are an open $O$ containing $P$ and a finite
family $(\phi_i)_{i \leq N}$ of isometric maps from $O$ to $K$ such that the
family of $(\phi_i (P))_{i \leq N}$ is disjoint. Thus: $m(\bigcup_i \phi_i (P)) = \sum_i m(\phi_i(P))
= N m(P)\leq 1$ using the invariance property. Hence, $m(P) \leq \frac{1}{N}$
for all non zero $N$, and $m(P) = 0$.

\begin{flushright}
$\blacksquare$
\end{flushright}

\section{Some consequences and applications of our results}

In this section, we show some applications of our existence and uniqueness theorems. These ones mainly deal with constructions of new measures in general enough cases to be out of
reach of usual methods. Moreover, some usual measures are recognized as given by our existence theorem. We study in a first time a very common case:
$K$ is supposed to be a compact subset of $\mathbf{R}^d$. We will see that most of the time, our measures are just probabilistic versions of the Lebesgue measure. That is why we state that they generalise it in a natural way.

\subsection{The case of finite dimensional numeric spaces}

\begin{theorem}
Let $K$ be a compact subset of $\mathbf{R}^d$ with a non-empty interior. It is assumed to be endowed with a translation invariant distance defining its usual topology. Then $\frac{\lambda}{\lambda (K)}$ is the only strongly isometrically compatible probability on $K$, where $\lambda$ denotes the Lebesgue measure. As a consequence, $\mu_K^L = \frac{\lambda}{\lambda (K)}$ for any limit notion $L$.
\end{theorem}

\vspace{2 mm}

\underline{\textbf{Demonstration:}} Any open set of $\mathbf{R}^d$ is homogeneous. Indeed, for any $x$ and $y$ in an open set, there is a $\varepsilon > 0$ small enough for balls centred in $x$ and $y$ with radius $\varepsilon$ to be included in the open, and so to be isometric. Moreover, if $x \in K$ and $y$ is in the interior of $K$, we easily see that : $x \hookrightarrow y$. Hence, our uniqueness theorem can apply. Denote $\mu$ this unique probability.

Now, $\mu$ is translation invariant in the interior of $K$, which implies that it equals $l\lambda$ for some $l > 0$ (use the uniqueness theorem for the Lebesgue measure in open sets) in $\overset{\circ }{K}$. Then, using the same arguments as in the proof of our strengthened uniqueness theorem, but only with translations, we get the conclusion.
\begin{flushright}
$\blacksquare$
\end{flushright}

The previous proof can be a bit disturbing. In fact, we are obliged to make this demonstrations with translations and not with isometric maps. Indeed, on some subsets of $\mathbf{R}^d$, partial isometries can have very strange behaviours. As a consequence, directly showing that $\lambda$ is a strongly isometrically compatible measure can be really difficult.

For instance, look at the following isometry. Let $K = \lbrace (x; \vert x\vert), x \in [-1;1]\rbrace$, $K' = [-1;1]\times \lbrace 0\rbrace$, and $f : K' \rightarrow K$ defined by : $f(x;0) = (x;\vert x\vert)$. Both $K$ and $K'$ are assumed to be endowed with the distance $d((x;y);(x';y')) = \max (\vert x-x'\vert; \vert y - y'\vert).$ We can show that $f$ is an isometry, but it is not the restriction of a linear map of $\mathbf{R}^2$.

\subsection{Application to Hausdorff measure}

Finding a general link between our probabilities and Hausdorff measures of compact sets of finite dimensional numeric spaces is a difficult question.

However, considering a compact set $K$ satisfying the conditions of one of our uniqueness theorems, and assuming that the Hausdorff measure $dh^{\alpha}$ of this compact is finite and non-zero, our theorem shows that $\mu_K$ and $\frac{dh^{\alpha}}{%
dh^{\alpha}(K)}$ equal. To show that, it is enough to prove that $%
dh^{\alpha} $ is strongly isometrically compatible, which raises no problem. More generally, that proves that under such hypothesis, the normalized Hausdorff measure is the unique strongly isometrically compatible probability on $K$.

Hence, our measure is often linked with Hausdorff measure. One can see here an interesting point of our construction: building a Hausdorff measure can be difficult, since the Hausdorff dimension is not necessarily a real number. It can be some kind of function, as can be seen in \cite{contrex}. Consequently, defining a probability measure using the Hausdorff process can be very difficult, and showing that it is always possible appears to be out of range. Our construction shows that a metrically compatible
measure exists when the space is metric compact. This result is relevant in this context.

\subsection{Link with the Haar measure}

Now, we exhibit the link between isometrically
compatible probabilities and the normalized Haar one on compact topological groups. We
will see that, once a good distance defining the topology of the compact
group is chosen, any such probability equals the Haar one, under
the assumption that the compact group is metrizable.

\begin{lemma}
Let $G$ be a topological metrizable compact group. There is a left and right
invariant distance on $G$ defining its topology.
\end{lemma}

\vspace{2 mm}

\underline{\textbf{Demonstration:}} Let $d$ be a distance defining the
topology of $G$. We put for $x$ and $y$ in $G$: $D(x;y) = \sup_{g\in
G ,g^{\prime }\in G} d(gxg^{\prime },gyg^{\prime }).$ The definition of $D$
makes sense by continuity of $(g;g^{\prime }) \mapsto d(gxg^{\prime
},gyg^{\prime })$ on a compact set. Moreover, $D$ is a distance (check the three axioms, there is no difficulty). And since $G$ is a group, one shows easily
that $D$ is invariant by left and right translations. Then, we have to show
that $D$ defines the same topology as $d$. First of all, it is obvious that $%
d(x;y) \leq D(x;y)$ for all $x$ and $y$. Let us show that if any sequence $(x_n)$
approaches some $x$ in $G$, that is $d(x_n;x)$ approaches $0$, then $D(x_n;x)$
approaches $0$ too. But for all $n$, $D(x_n;x)$ is a supremum of a
continuous function on a compact set. Hence for all $n$, there are $g_n$ and 
$g^{\prime }_n$ such that $D(x_n;x) = d(g_nx_ng^{\prime }_n,g_nxg^{\prime }_n)$. Let then $L$ be a limit point of the sequence $(D(x_n;x))$. There is
a strictly increasing sequence $\varphi$ of natural numbers such that $%
(D(x_{\varphi (n)};x))$ approaches $L$. But $(g_{\varphi (n)})$ and $%
(g^{\prime }_{\varphi (n)})$ being sequences in a compact, there is a
strictly increasing sequence $\psi$ of natural numbers such that $%
(g_{\varphi (\psi (n))})$ and $(g^{\prime }_{\varphi (\psi (n))})$ approach,
say, respectively, $g$ and $g^{\prime }$. Then $(D(x_n;x))$ approaches $d(gxg^{\prime };gxg^{\prime }) = 0$ by continuity of $d$.
Hence $L = 0 $ and $(D(x_n;x))_n$ has only one limit point: $0$. Thus $D(x_n;x)$ approaches zero, and $D$
defines the same topology as $d$.

\begin{flushright}
$\blacksquare$
\end{flushright}

Consequently :

\begin{theorem}
Let $G$ be a metrizable topological compact group and $D$ a distance on $G$
defining its topology and being invariant under left and right translations.
Then any isometrically compatible probability of $(G;D)$ equals the Haar probability of $G$.
\end{theorem}

\vspace{2 mm}

\underline{\textbf{Demonstration:}} Let $m$ be an isometrically compatible
probability of $G$. It is invariant under partial isometries on opens, hence
under left and right translations on opens. Let then $B$ be a Borel set of $%
G $ and $g$ an element of $G$. We have for all open $O$ containing $B$: $%
m(O) = m(gO) = m(Og).$ Since probabilities defined on Borel sets are
regular, the fact that $x \mapsto gx$ and $x \mapsto xg$ are bi-continuous
implies that: $m(B) = m(gB) = m(Bg)$ using infima. $m$ hence
satisfies the characteristic properties of the probabilistic Haar measure on 
$G$. Since there is only one such measure, our measures equal.
\begin{flushright}
$\blacksquare$
\end{flushright}

Note that the usual proof of existence of the Haar measure is not constructive. Our result is interesting since it gives informations about the effective calculus of this measure when the difficult question of the manipulation of the limit notion can be simplified into a classical limit calculus.

\subsection{A patching measures theorem and application to locally compact homogeneous metric spaces}

Now, we show in this section that homogeneous locally compact metric space can be endowed with a Borel locally finite measure wich is invariant under partial isometries definde on opens. The strategy for the demonstration is a patching process of the measures given by the Geometry and Probability theorem using the following result. This result has its own interest and can probably be used in other contexts.

\begin{theorem} 
Let $E$ be a set endowed with a $\sigma$-algebra $M$ such that there is a family $(B_i)_{i \in I}$ of measurable sets covering $E$, and such that for all $i\in I$ there is a non-negative measure $m_i$ on $B_i$ defined on the induced $\sigma$-algebra. There is then a measure $\mu$ on $E$ defined on $M$, canonically chosen, such as: 
$$\mu (B) = \sup \sum_{i} m_i(A_i),$$
where the supremum is taken over disjoint families $(A_i)_{i \in I}$ of measurable sets contained in $B$, with for all $i$, $A_i \subset B_i$, and such that the set of $i$ satisfying $A_i \neq \emptyset$ is countable. It is called the supremum measure of the $(m_i)_{i \in I}$.

Moreover, assuming that for all $i$ and $j$ in $I$, and all measurable set $B$ of $M$ contained in $B_i\cap B_j$, $m_i(B) = m_j(B)$ then for all $i\in I$, and all measurable set $B$ contained in $B_i$, $\mu (B) = m_i(B)$. This last measure is called the canonical patching measure of the family $(m_i)_{i \in I}$.
\end{theorem}

\vspace{2 mm}

The fact that no cardinal assumption has been made about $I$ has to be pointed out.

\vspace{2 mm}

\underline{\textbf{Demonstration:}} Let us use the definition of the theorem: for all measurable set $B$ of $E$, 
$$\mu (B) = \sup \sum_{i} m_i(A_i),$$
where the supremum is taken over disjoint families $(A_i)_{i \in I}$ of measurable sets contained in $B$, with for all $i$, $A_i \subset B_i$, and such that the set of $i$ satisfying $A_i \neq \emptyset$ is countable.

Let us show that $\mu$ is a measure. Let $(C_j)_{j \in \mathbf{N}}$ be a countable family of disjoint measurable sets, and $C$ its union. We have to prove that $\mu(C) = \sum_{j} \mu (C_j)$.

Choose a $\varepsilon >0$. There is a family $(A_i)_{i \in I}$ of measurable sets contained in $C$, with for all $i$, $A_i \subset B_i$, such that the set of $i$ satisfying $A_i \neq \emptyset$ is countable, and such as: $\mu(C) \geq \sum_{i} m_i(A_i) \geq \mu(C) - \varepsilon.$ Then: $\mu(C) - \varepsilon \leq \sum_{i} \sum_{j} m_i(A_i \cap C_j) =  \sum_{j} \sum_{i} m_i(A_i \cap C_j) \leq \sum_{j} \mu(C_j)$. Indeed, one can swap the sum symbols in series whose terms are non-negative, and the $(A_i \cap C_j)_{i \in I}$ are disjoint measurable sets contained in $C_j$, with for all $i$, $A_i \cap C_j \subset B_i$, and such that the set of the $i$ satisfying $A_i\cap C_j \neq \emptyset$ is countable. This being true for all $\varepsilon >0$, we get an inequality.

Then, a $\varepsilon >0$ still being chosen, for all $j$, there is a disjoint family $(C_{i,j})_{i \in I}$ of measurable contained in $C_j$, with for all $i$, $C_{i,j} \subset B_i$, such that the set of $i$ satisfying $C_{i,j} \neq \emptyset$ is countable, and such as: $\mu (C_j) \leq \sum_i m_i(C_{i,j}) + \frac{\varepsilon}{2^j}.$ Then: $ \sum_j \mu (C_j) \leq 2\varepsilon + \sum_j \sum_i m_i(C_{i,j}) = 2\varepsilon + \sum_i \sum_j m_i(C_{i,j}) = 2\varepsilon + \sum_i m_i(\coprod_j C_{i,j}).$ Indeed, the $(C_{i,j})_{j \in \mathbf{N}}$ are disjoint and contained in $B_i$ for all $i$. Then: $ \sum_j \mu (C_j) \leq 2\varepsilon + \sum_i m_i(\coprod_j C_{i,j}) \leq \mu (C) + 2\varepsilon,$ since the $(\coprod_j C_{i,j})_{i \in I}$ are disjoint measurable sets contained in $C$, such that for all $i$, $\coprod_j C_{i,j} \subset B_i$, and such that the set of $i$ satisfying $\coprod_j C_{i,j} \neq \emptyset$ is countable. Indeed, $\coprod_j C_{i,j} \neq \emptyset \Leftrightarrow \exists j \in \mathbf{N} / C_{i,j} \neq \emptyset$, hence the set of $i$ such that $\coprod_j C_{i,j} \neq \emptyset$ is a countable union of countable sets. This being true for all $\varepsilon >0$, we get the reverse inequality.

This shows that $\mu$ is a measure. Let us prove now that it prolongs the $(m_j)$ under the last assumption: choose a $j$ in $I$, and $B$ a measurable set of $B_j$. First of all, $\mu(B) \geq m_j(B)$ (look at the definition). Then, for all $\varepsilon >0$, there is a disjoint family $A_i$ of measurable sets contained in $B$, such that for all $i$, $A_i \subset B_i$, and such that the set of $i$ satisfying $A_i \neq \emptyset$ is countable, with:  $\mu(B) \leq \sum_{i} m_i(A_i) + \varepsilon.$ But all the $A_i$ are contained in $B_j$, and with the theorem assumption: $m_i(A_i) = m_j(A_i)$. Hence, with $\sigma$-additivity, $$\mu(B) \leq m_j(\coprod_{i} A_i) + \varepsilon \leq m_j(B) + \varepsilon.$$ This being true for all $\varepsilon >0$, we get the conclusion.

Finally, this construction is canonical, since it uses no arbitrary choice.\begin{flushright}
$\blacksquare$
\end{flushright}

Note that when $I$ is countable, that patching measure is unique: to prove it, just build a countable measurable partition of $E$ whose elements are contained in the $B_i$. By instance, choose $P_n = B_n\setminus (\bigcup_{k\leq n-1}B_k).$ One can then show the uniqueness proving that for a measurable set $B$, the only possible measure $\mu(B)$ is: $\mu(B) = \sum_n m_n(B \cap P_n).$

At last, note that when $I$ is not countable, many measures can be compatible with restrictions. Indeed, chose $E = \mathbf{R}$ and $B_i = \lbrace i\rbrace$ for all $i \in \mathbf{R}$, and $m_i = 0$. Then the Lebesgue measure and the zero measure are compatibles with this system.

\vspace{2 mm}

Now let us show the main result of this paragraph.

\begin{theorem} 
Let $K$ be a locally compact homogeneous metric space. There is a Borel non-null locally finite non-negative measure $P$ on $K$ which is invariant under partial isometries defined on open sets of $K$.
\end{theorem}

\vspace{2 mm}

\underline{\textbf{Demonstration:}} Choose for each $x \in K$ a ball $O_x$ centred in $x$ and whose closure is compact. Our main theorem gives a probability $P'_x$ on each $\overline{O_x}$, and $(O_x)_x$ covers $K$. Now, choose a $x_0\in K$. For each $x \in K$, there are a real $r_{x_0;x} > 0$ and an isometry $T$ defined on $B(x_0;r_{x_0;x})$ such that $x \in T(B(x_0;r_{x_0;x}))$, and $P'_x$ and $T^{-1}_{*}P'_{x_0}$ are isometrically invariant measures defined on the homogeneous pre-compact space $T(B(x_0;r_{x_0;x}))$. Hence, there is a $l_x$ such that $P'_x = l_x T^{-1}_{*}P'_{x_0}$ on the Borel sets of $T(B(x_0;r_{x_0;x}))$ and $l_x > 0$ since $P'_x$ is non null (otherwise, it would not be a probability). Denote $P_x = \frac{1}{l_x} P'_x$.

Now, let us show that the system $(P_x, \overline{O_x})$ is patching. Let $x$ and $y$ be any elements of $K$, and $F = \overline{O_x} \cap \overline{O_y}$. If $F = \emptyset$, the patching condition is true. Otherwise, $P_x$ and $P_y$ are isometrically invariant measures defined on the homogeneous pre-compact space $F$. Hence there is a $l'$ such that $P_x = l' P_y$ on $U$, and one can show that $l' = 1$. Indeed, choose $z \in U$. A well chosen ball $B$ centred in $z$ is isometric to a Borel subset of another centred in $x$, and to a Borel subset of another centred in $y$, and finally to a Borel subset of another centred in $x_0$. Choosing this ball small enough, one sees that $l' = 1$ if $P_y (B) \neq 0$. Now, if for any ball $B$ of $U$, $P_y (B) = 0$, with compactness, $P_y(U) = 0$, and one can choose $l' = 1$ too.

Thus the system $(P_x, \overline{O_x})$ is patching. Let $P$ be the patched measure. It is obviously non-null and locally finite. Let $T$ be a partial isometry between open sets of $K$ defined on, say, $O$. $P$ and $T_{*}P$ equal on each $\overline{O_x}\cap T^{-1}(\overline{O_y})$ since they are isometrically invariant measures defined on a homogeneous pre-compact space equalling on a non-null Borel set. To show this last point, just make the same kind of proof as before. Now, if $O$ can be written as a countable union of non-null $\overline{O_x}\cap T^{-1}(\overline{O_y})$, with $\sigma$-additivity, $P(O) = P(T(O))$, and if not, $P(O) =+\infty = P(T(O))$. Consequently $P$ is invariant under partial isometries between open set.\begin{flushright}
$\blacksquare$
\end{flushright}

Moreover, our uniqueness theorem enables us to state :

\vspace{2 mm}

\begin{theorem} 
Let $K$ be a locally compact homogeneous metric space. Let $P$ be a Borel non-null locally finite non-negative measure $P$ on $K$ which is invariant under partial isometries defined on open sets of $K$, and $Q$ be another. Then there is a $l>0$ such that $Q = l P$. 
\end{theorem}

Indeed, one easily show that such a $l$ exists on each ball whose closure is compact, and studying their intersections, one sees that all these $l$ must equal.

\subsection{A version of our theorem for pseudo-metric compact spaces}

We will see in this paragraph that our main theorem can be generalized to compact pseudo-metric spaces.

\vspace{2 mm}

\begin{theorem} 
Let $K$ be a compact pseudo-metric space. There is a probability measure $P$ defined on the Borel $\sigma$-algebra of $K$ such that partially isometric opens have the same probability.
\end{theorem}

In fact, one can show that the whole theorem can be stated in this new framework. But the previous version will be sufficient for the future application we aim at.

\vspace{2 mm}

\underline{\textbf{Demonstration:}} Let $K'$ be the quotient space of $K$ by the kernel of the pseudo-distance of $K$. $K'$ is 
compact metric, hence it admits a metrically compatible probability $P'$. Denote $\pi$ the canonical map from $K$ to $K'$. It is an 
isometry, and it induces a $\sigma$-algebras morphism between the Borel sets of $K'$ and the saturated Borel sets of $K$. Indeed, if 
$A$ and $B$ are saturated subsets of $K'$ (that is $A = \pi^{-1}(\pi (A))$ and $B = \pi^{-1}(\pi (B))$), $\pi (A \cup B) = \pi (A) 
\cup \pi (B)$, and $\pi (A \cap B) = \pi (A) \cap \pi (B)$. For the last point, the only difficulty is to show that $\pi (A \cap B) 
\supset \pi (A) \cap \pi (B)$. Let $y \in \pi (A) \cap \pi (B)$. There is $x \in A$ and $x' \in B$ such that $y = \pi (x) = \pi (x')
$. But $A$ and $B$ are saturated. Since $y \in \pi (B)$ and $y = \pi (x)$, $x \in B$, thus $x \in A \cap B$, and $y \in \pi (A \cap 
B)$.
Moreover, the previous properties of $\sigma$-algebras morphism can be stated for unions and intersections of any cardinal.

Hence, the following formula : $$P(A) = P'(\pi (A))$$ defines a probability measure on the $\sigma$-algebra of saturated Borel subsets of $K$, and this probability is such that partially isometric saturated opens have the same probability. Indeed, if $f$ is a partial isometry of $K$, defining $\bar{f}$ by $\bar{f} ([x]) = [f(x)]$, one gets a well defined isometry of $K'$. The demonstration is straightforward once this is noticed.

Finally, the $\sigma$-algebra of saturated Borel subsets of $K$ and the $\sigma$-algebra of Borel subsets of $K$ are the same. Indeed, both are generated by open balls of $K$ (since $K$ has a countable dense subset). To show that, just observe that any open ball is saturated thanks to its definition by the pseudo-distance.\begin{flushright}
$\blacksquare$
\end{flushright}

\subsection{Probability on the space of probabilities defined
on the Borel sets of a non-empty Hausdorff compact space}

In this section, we construct a probability on the space of probabilities
defined on the Borel sets of a non-empty Hausdorff compact space. This
measure is invariant by numerous maps.

Denote $K$ a non-empty Hausdorff compact space, $C$ the set of continuous
maps from $K$ to $\mathbf{R}$, and $M$ the set of Radon measures on $K$. $C$
is endowed with the norm of the supremum. Following the Riesz Theorem, $%
M $ equals the topological dual of $C$. Finally, one denotes $P$ the set of
probabilities defined on Borel sets of $K$.

$M$, as dual of $C$, is endowed with the dual norm: 
\begin{equation*}
\Vert \mu \Vert = \sup_{\Vert f\Vert_{\infty} \leq 1} \vert \int_K f d \mu
\vert. 
\end{equation*}
It is obvious that $P$ is included in the unit ball of $M$ associated with
this norm.

Moreover, $M$ is endowed with the weak-$*$ topology since it is the
dual of $C$. And $C$ being separable, the unit ball of $M$ is metrizable for the
weak-$*$ topology (see \cite{Brézis}[page 48]). If $(f_n)$ is a denumerable
dense subset of the unit ball of $C$, the following distance defines the
weak-$*$ topology: $d^*(\lambda , \mu ) = \sum_n \frac{1}{2^n}\vert \int_K
f_nd \lambda - \int_K f_nd \mu\vert$.

Finally, $P$ is closed in $M$ for the weak-$*$ topology. Indeed, $P =
\lbrace \mu \in M / \int_K d\mu = 1 \mbox{ and } \forall f \in C, \int_K
fd\mu \geq 0\rbrace$. Thus $%
P$ is a weak-$*$ compact. Consequently, there is a metrically compatible
probability $m$ constant on partially isometric opens for the chosen weak-$*$
distance of $P$.

Let us try to construct some partial weak-$*$ isometries of $P$, in order to
see that our result is useful.

First of all, it is obvious that $d^*$ is translation-invariant in the
following sense: if $B$ is a weak-$*$ Borel set of $P$ such that for some
measure $p$, $B+p$ is included in $P$, then for all measures $\lambda$ and $%
\mu$ of $B$, $d^*(\lambda ; \mu) = d^*(\lambda +p ; \mu +p)$. Necessarily, $%
p(K)=0$ for such measures, hence they are not non-negative if non-zero.

It follows that $m$ is invariant under translations between opens, if well
defined.
\vspace{2 mm}

\begin{theorem}
Let $K$ be a compact topological space, and $M(K)$ be the set of probabilities defined on Borel subsets of $K$, endowed with the weak-$*$ distance $d^*(\lambda , \mu ) = \sum_n \frac{1}{2^n}\vert \int_K
f_nd \lambda - \int_K f_nd \mu\vert$ where $(f_n)$ is a denumerable
dense subset of the $\Vert .\Vert_{\infty}$-unit ball of $C(K;\mathbf{R})$.
There is a probability $m$ on $M(K)$ which is invariant under the isometries of $M(K)$ defined on open sets.
\end{theorem}

\vspace{2 mm}

Here is another interesting construction, when $K$ is assumed to be metric.
Denote $d$ a distance defining the topology of $K$. One can
consider its isometry group $G$, with the supremum distance. With the
Ascoli Theorem, one sees that this group is compact. So, it possesses a Haar
measure $dh$. Let us define a new distance on $P$ by:

\begin{equation*}
D(\lambda ; \mu) = \int_G d^*(\phi_{*} (\lambda) ; \phi_{*} (\mu)) dh(\phi), 
\end{equation*}
where the notations $\phi_{*} (\lambda)$ and $\phi_{*} (\mu)$ denote the
push-forward measures.

$D$ is well defined, since $\phi \mapsto d^*(\phi_{*} (\lambda) ; \phi_{*} (\mu))$
is integrable. Indeed, it is bounded, and we will show that it is measurable.

Let $\phi$ be an isometry of $G$ and $B$ be a closed set of $K$. Let $(\phi
_n)$ be a sequence of $G$ approaching $\phi$. It is obvious that $%
(d_H(\phi^{-1} (B);\phi _n^{-1}(B)))_n$ approaches zero ($d_H$ is the
Hausdorff distance associated to $d$). Now, one knows that any finite positive measure defined on Borel sets
is upper-semi-continuous with respect to closed sets endowed with the
Hausdorff distance (it is a consequence of regularity). Hence, here: $%
\limsup_{n \rightarrow +\infty} \lambda (\phi _n^{-1}(B)) \leq \lambda (\phi
^{-1}(B))$. So $\phi \mapsto \int_K \textbf{1}_B d(\phi_{*} (\lambda))$ is measurable for any closed set $B$ and any probability $\lambda$. But, for any Borel set $B'$, 
$$\lambda (\phi^{-1} (B')) = \sup_{F \subset B', F \mbox{ closed}} \lambda (\phi^{-1} (F)),$$
and this supremum can be taken among a countable set. 
Hence, $\phi \mapsto \int_K \textbf{1}_{B'} d(\phi_{*} (\lambda))$ is measurable for any Borel set $B'$. It follows that the same property holds for
integrals of simple non-negative functions, hence with a supremum, for
integrals of non-negative (integrable) functions. Consequently, for any non-negative
integrable function $f$ and any probability $\lambda$ on $K$, $\phi
\mapsto \int_K fd(\phi_{*} (\lambda))$ is Borel measurable. It is thus still true for all
integrable functions. Then, it is enough to look at the definition of $d^*$
to see that $\phi \mapsto d^*(\phi_{*} (\lambda) ; \phi_{*} (\mu))$ is measurable
for any probabilities $\lambda$ and $\mu$.

Then, $D$ is a distance (the proof is easy). Finally, for all isometries $%
\varphi$ of $K$, and all measures $\lambda$ and $\mu$ of $P$, $D(\lambda ;
\mu) = D(\varphi_{*} (\lambda) ; \varphi_{*} (\mu))$. To show that, just use the
characteristic property of $dh$.

Hence, the map $\lambda \mapsto \varphi_{*} (\lambda)$ is a $D$-isometry of $P$.

Let us show at last that $P$ is compact for $D$.

Let $(\lambda_n)$ be a sequence of $P$. Since $P$ is compact for $d^*$, this
sequence is such that there is a strictly increasing sequence of natural
numbers $\psi$ with $(\lambda_{\psi (n)})$ approaching some $\lambda$ of $P$%
. Looking at the definition of the weak-$*$ topology, we see that it implies
that for any isometry $\phi$, $(\phi_{*} (\lambda_{\psi (n)}))$ approaches $%
\phi_{*} (\lambda )$. Thus $(d^*(\phi_{*} (\lambda) ; \phi_{*} (\lambda_{\psi (n)})))_n$
simply approaches zero and is bounded. With Lebesgue theorem, $%
(D(\lambda_{\psi (n)};\lambda))_n$ approaches zero.

It is then enough to consider a metrically compatible probability $m$
given by the main theorem for the compact $P$ endowed with $D$. It is
invariant by partial maps $\lambda \mapsto \varphi_{*} (\lambda)$ where $\varphi$ is
an isometry of $K$. Moreover, one can easily check that this measure is
invariant under well defined partial translations of $P$.
\vspace{2 mm}

\begin{theorem}
Let $K$ be a compact metric space, and $M(K)$ be the set of probabilities defined on Borel subsets of $K$, endowed with the weak-$*$ topology.
There is a probability $m$ on $M(K)$ which is invariant under the partial maps $\lambda \mapsto \varphi_{*} (\lambda)$ defined on weak-$*$ opens where $\varphi$ is an isometry of $K$.
\end{theorem}

\vspace{2 mm}

The previous definition of a relevant probability on the space of probabilities on a compact set can appear to be very abstract. Actually, we can give us a concrete representation of what it is.

Imagine a technical process consisting in a machine throwing exactly the same mass of dust in Petri dishes, but with a random distribution. Each distribution of dust can be modelled by a probability on the compact Petri dish, and we need a probability on the space of dust distributions to model this experiment. Our measure can be used.

\subsection{Probability measure on some spaces of subsets
of a compact metric space}

We define here a notion of metrically compatible probability on the space
of closed non-empty subsets of a compact metric space $K$. Denote $[K]$ this
set. It is a compact for the Hausdorff distance. It is then enough to apply
our main theorem. Moreover, one can make the same kind of construction with
the opens of a compact metric space different from the whole space. Just endow their set $T\setminus \lbrace K\rbrace$ with the
following distance: $d(O,O^{\prime }) = d_H(\complement_K^O;\complement_K^{O^{\prime }})$,
where $d_H$ is the Hausdorff distance. Indeed, the map $O \mapsto \complement_K^O$ is
a bijection between $T$ and $[K]$.

As above, that construction is really interesting only if one can exhibit
many partial isometries. Some follow.

Let $\phi$ be a partial isometry between two opens $O$ and $O^{\prime }$ of $%
K$, and $U$ an open included in $O$. Denote $[U]$ the set of closed of $K$
contained in $U$, and $[\phi (U)]$ its image by $F \mapsto \phi (F)$, which
is defined on closed sets of $O$. That map is an isometry for the Hausdorff
distance (the demonstration is easy), and those sets are open for that
distance. Indeed, for any closed $F$ included in, say, $U$, $\alpha =
d(F;\complement_K^U)$ is non-zero because we are in a compact. Then, if $d_H(F,
G)<\alpha$ where $G$ is any closed of $K$, one has: $G \subseteq U$ (check this using the definition of the Hausdorff distance $d_H$).

Thus, the sets $[U]$ and $[\phi (U)]$ have the same probability.

\vspace{2 mm}

Now, we present a fact which can be felt quite disturbing, but which, well understood, enables to see the real importance of the previous construction. This result has been pointed out to me by M. Drouin. Choose $K = [0;1]$. One will show that $\mu_{[K]} = \delta_{K}$, the Dirac measure on the full subset $K$. Let $\alpha >0$ and $c\in [0,1]$, and try to calculate $\mu_{[K]}([K\setminus \overline{O}])$ where $O = B(c,\alpha)$. $K\setminus \overline{O}$ is open in $[K]$, hence: $\mu_{[K]}([K\setminus \overline{O}])\leq
L_{\varepsilon\rightarrow 0} \frac{N([K\setminus \overline{O}];\varepsilon)}{N([K];\varepsilon)}$.

But $K$ can be covered with at least $E(1/\varepsilon)+1$ open balls of radius $\varepsilon /2$ (where $E(x)$ denote the tallest integer lower than $x$). Hence $[K]$ can be covered with at least $2^{E(1/\varepsilon)+1}-1$ balls of radius $\varepsilon /2$. Indeed, if $(B(x_i,\varepsilon /2))_{i\leq N}$ is a covering of $K$, $([\bigcup_{i\in P} B(x_i,\varepsilon /2)]_{P \subset [[0;N]], P \neq \emptyset})$ is a covering of $[K]$, and the $[\bigcup_{i\in P} B(x_i,\varepsilon /2)] = B_{[K]}(\lbrace x_i, i \in P\rbrace , \varepsilon /2)$ are balls in $[K]$. In order to see that, choose $A$ a subset of $K$, and take for $X$ the largest subset of $\lbrace x_i, i \leq N\rbrace$ such that $d(x_i, X) < \varepsilon /2$ for all $i$. Since $(B(x_i,\varepsilon /2))_{i\leq N}$ is a covering, $A \subset \bigcup_{x \in X} B(x;\varepsilon /2)$ and $X \subset \overline{A}^{\varepsilon /2}$, hence $A$ belongs to $B_{[K]}(\lbrace x, x \in X\rbrace , \varepsilon /2)$.

Now, let us show that $[K]$ cannot be covered with less that $2^{E(1/\varepsilon)+1}-1$ balls of radius $\varepsilon /2$. Let $x_i = i\varepsilon$ for $i \leq E(1/\varepsilon)$ of $K$, and denote $L$ their set. Any two distinct non-empty subsets of $L$ have a distance in $[K]$ at least equal to $\varepsilon$. Now, let $(B_i)_i$ be any covering of $[K]$ by open balls of radius $\varepsilon /2$. Any non-empty $P \subset L$ belongs to some $B_i$ since $(B_i)_i$ is a covering of $[K]$, and if $P \subset L$ and $ Q \subset L$ belong to the same $B_i$, $P=Q$. Indeed, if, for instance, there is a $x \in Q \setminus P$, $d_H(Q;P)\geq \varepsilon$, whence for all $F$ and $F'$ in $B_i$, $d_H(F;F') < \varepsilon $. So $(B_i)_i$ has at least as many balls as the cardinal of the family of non-empty subsets of $L$, that is $2^{E(1/\varepsilon)+1}-1$.

Therefore $N([K];\varepsilon) = 2^{E(1/\varepsilon)+1}-1$. And with the same kind of reasoning, one 
has: $N([K\setminus \overline{O}];\varepsilon)= 2^{E(\frac{1-\alpha}{\varepsilon})+1}-1$. Hence: 

$$L_{\varepsilon\rightarrow 0} \frac{N([K\setminus \overline{O}];\varepsilon)}{N([K];\varepsilon)} \leq
L_{\varepsilon\rightarrow 0}\frac{2^{E(\frac{1-\alpha}{\varepsilon})+1}-1}{2^{E(1/\varepsilon)+1}-1} =
L_{\varepsilon\rightarrow 0} 2^{E(\frac{1-\alpha}{\varepsilon})-E(1/\varepsilon)}=0,$$

and $\mu_{[K]}([K\setminus \overline{O}])=0$.

Consequently, $\mu_{[K]}(\bigcup_{c \in [0,1], \alpha > 0}[K\setminus \overline{B(c,\alpha)}])=0$ ($K$ is separable). But any closed of $K$ different of $K$ is contained in $\bigcup_{c \in [0,1], \alpha > 0}[K\setminus \overline{B(c,\alpha)}]$. Hence $\mu_{[K]} ([K]\setminus \lbrace K\rbrace) =0$, and $\mu_{[K]}$ is concentrated on $K$, which ends the demonstration.

\vspace{2 mm}

One sees then that in the case of a sufficiently homogeneous space $K$, the measure $\mu_{[K]}$ is concentrated on the whole space (as a matter of fact, the previous reasoning holds in many cases, provided the space satisfies good hypothesis). One can see on examples that when the space is not homogeneous, the support of the measure $\mu_{[K]}$ gives informations about the homogeneity of $K$. More precisely one can consider the spaces of metrically compatible probabilities on $[K]$. For instance, denote $K_1$ and $K_2$ two compact metric spaces such that no $1$-coercive nor $1$-lipschitz map exist between no-empty opens of $K_1$ and $K_2$ (see the M. Drouin's counter-example of the beginning of our paper for example. There is some work to show that such maps do not exist.) Then $\lambda_c (A) = c\mu_{[K_1]}(A\cap [K_1])+(1-c)\mu_{[K_2]}(A\cap [K_2])$ defines a metrically compatible probability on $[K_1 \amalg K_2]$ for $c \in [0;1]$, and it is not a Dirac measure when $c \in (0;1)$. The space $\lbrace \lambda_c, c \in [0;1]\rbrace$ contains informations about homogeneity properties of $K_1 \amalg K_2$ : its extreme points are Dirac measures on subsets of $K_1 \amalg K_2$ such that both are homogeneous, and such that no $1$-coercive nor $1$-lipschitz partial map between opens of these exist.

Hence, given a compact metric space $K$, one can consider the spaces of invariant probabilities of $K$ (known to be non empty), either metrically, strongly isometrically or isometrically compatible, and then the spaces of invariant probabilities of $[K]$ (non empty too), and so on. These convex spaces of probabilities contain informations about the metric topology of $K$, and some theorems giving relations between these spaces and the metric topology of $K$, peculiarly its homogeneity characteristics must exist. We do not present these theorems, but just suggest that a theory dealing with those questions can be built.

\subsection{Application of the previous construction to phylogeny}

The previous construction of probability measures on the space of closed subsets of a compact metric space can be applied to build a model in phylogeny. Indeed, the evolution of species from a common ancestor is described by a tree. Practically, one can imagine a maximal tree modelling the most complex evolution from the common ancestor (figured by the root of the tree), but in reality, the observed evolution is only a sub-tree of this maximal tree. Actually, the real evolution tree can be seen as a random tree chosen among the sub-trees of the maximal possible tree. The probability law describing the choice of this sub-tree has to possess invariance properties. Indeed, if two families of sub-trees are the same up to an exchange of the mutations loci, and if these mutations have the same probabilities, the two families of sub-trees must have the same probabilities.

Here, we give a very simplified version of a solution to this problem. In fact, a much better solution can be built with the same ideas and the same mathematical tools : our main theorem and a trick enabling to see sub-trees as sub-spaces of the space of branches. A better solution would have to incorporate in the model the probability of each mutation and the delay between two of them.

\vspace{2 mm}

Let us describe this model. Let $K = \lbrace 0 ; 1\rbrace^{\mathbf{N}}$ be the set of branches of our maximal tree, taken binary for 
simplicity, and chosen such that every mutation has the same probability. This tree is : 
$$T = \lbrace \mbox{finite sequences of } \lbrace 0;1\rbrace\rbrace,$$
endowed with the tree-order relation $ (a_n)_{n\leq N} \leq (b_n)_{n\leq P}$ if $N \leq P$ and $b_n = a_n$ for all $n \leq N$. This tree is usually metrised by $d(x;y) = e^{-h(x\wedge y)}$, where the height $h$ gives the number of 
element of a finite sequence, and $\wedge$ is the infimum for the order relation. We do not check the usual results about all these 
structures : it can be shown that $\wedge$ exists, and that $d$ is a distance. Moreover, $d$ can be defined for branches of $K$, and 
one can easily check that $d$ is a way to define the component-wise convergence in $K$.

Now, actual evolution trees will be taken as complete sub-trees of $T$, that is we make the assumption that there is no disappearance of specie. With this hypothesis, an evolution tree is a subset $S$ of $T$ such that $x \in S$ and $y\leq x$ implies $y \in S$ (sub-tree condition), and such that $x\in S$ implies that there is a $y\in S$ such that $y \geq x$ and $y\neq x$ (completion condition). In this framework, it can be shown that evolution sub-trees of $T$ are characterised by the set of their branches into $K$, and that these sets are exactly the closed subsets of $K$.

Hence, our problem of randomly choosing an evolution sub-tree of $T$ can be transformed into a problem of a random choice of a closed subset of $K$, a compact metric space. This has already been done just before.

\vspace{2 mm}

Finally, we managed to build a probability law on evolution trees in our simplified model, and this law respects an invariance property under partial isometries for a distance linked to $d$ ( the $d$-Hausdorff distance). In fact, this invariance property insures us that the law has a physical meaning.

\subsection{Application to the Krylov-Bogolioubov problem}

In this section, we give a version of a general Krylov-Bogolioubov theorem in the case of a family of measurable maps acting on a compact metric space. I have read (from Mathoverflow) that an original version of that theorem, showed by Bogolioubov, deals with any amenable group of maps acting on a locally compact space. Unfortunately, this theorem, first written in Russian, has never been translated in English. I do not speak Russian, but I managed to apply our "Geometry and Probability Theorem" to a quite interesting related question. Moreover, this new theorem states a strong property : it gives some kind of equivalence for the existence of an invariant probability.

\vspace{2 mm}

\begin{lemma} Let $(K;d)$ be a compact metric space. There is a Borel probability measure on $K$ such that its support is $K$ itself.
\end{lemma}

\vspace{2 mm}

\underline{\textbf{Demonstration:}} Let $E = \lbrace x_n, n\in \mathbf{N}\rbrace$ be a countable dense subset of $K$, assuming $x_n \neq x_p$ if $n \neq p$. Denote $P = \sum_n \frac{1}{2^{n+1}}\delta_{x_n}$, where $\delta_{x_n}$ is the Dirac measure in $x_n$. $P$ is a Borel probability measure, and the largest open set whose measure is null is $\emptyset$.\begin{flushright}
$\blacksquare$
\end{flushright}

\begin{theorem}[Krylov-Bogolioubov Equivalence]
Let $(K;d)$ be a compact metric space.

\vspace{2 mm}

(1) Let $P$ be any probability defined on the Borel $\sigma$-algebra of $K$. Then the set $M = \lbrace f: K \rightarrow K \mbox{ Borel } / \forall A \in Bor(K), P(f^{-1}(A)) = P(A)\rbrace$ is a compact monoid for the composition law and an invariant under right-translations pseudo-distance $D$. Moreover, $f \mapsto d(f(x);f(y))$ is measurable for the Borel $\sigma$-algebra completed by null sets defined by $D$, for almost every $x$ and $y$ in $K$, and if $(f_n) \rightarrow f$ point-wisely, $(D(f_n;f)) \rightarrow 0$.

\vspace{2 mm}

(2) If $(S;D)$ is a pre-compact monoid of measurable transformations of $K$ for an invariant under right-translations pseudo-distance $D$ such that if $(f_n) \rightarrow f$ point-wisely, $(D(f_n;f)) \rightarrow 0$, there is a compact monoid $(M;D)$ of 
measurable transformations of $K$ for the composition law such that $D$ is an invariant under right-translations pseudo-distance, and $S \subset M$.

\vspace{2 mm}

(3) If $(M;D)$ is a compact monoid of measurable transformations of $K$ for the composition law such that $D$ is an invariant under right-translations pseudo-distance, there is a Borel probability measure $m$ on $M$ which is invariant under right-translations.

\vspace{2 mm}

(4) If $(M;D)$ is a compact pseudo-metric monoid of measurable transformations of $K$ endowed with a Borel probability measure $m$ which is invariant under right-translations for the composition law , and if $f \mapsto d(f(x);f(y))$ is measurable for the Borel $\sigma$-algebra completed by null sets defined by $D$ for almost every $x$ and $y$ in $K$, $K$ can be endowed with an invariant under $M$-action pseudo-distance $d'$ such that $(K;d')$ can be written as the union of a $d'$-$\sigma$-compact set $\Sigma$ with a null set for some Borel probability $\mu$ on $K$. Moreover, the topology defined by $d'$ on $\Sigma$ is included in the one defined by $d$ on $\Sigma$, and if $d'$ is a true distance, they equal.

\vspace{2 mm}

(5) Finally, assume that $(K;d)$ is a compact metric space such that for some $d'$ and some $d$-Borel probability $\mu$, $(K;d')$ is the union of a $\sigma$-compact pseudo-metric space with a $\mu$-null set, and that the $d'$-Borel sets are $d$-Borel sets. In addition, assume that $S$ is a set of measurable transformations of $K$ such that $d'$ is invariant under the action of the elements of $S$. Then there is a $d'$-Borel completed probability $P$ on $K$ which is invariant under every maps of $S$.

\vspace{2 mm}

Peculiarly, a simpler partial reverse of (1) follows :

(2') Let $(K;d)$ be a compact metric space. Assume that $(M;D)$ is a compact monoid of measurable maps of $K$ for the composition law and that $D$ is an invariant under right-translations distance such that $f \mapsto d(f(x);f(y))$ is measurable for the Borel $\sigma$-algebra completed by null sets defined by $D$ for almost every $x$ and $y$ in $K$. Then there is a Borel probability $P$ on $K$ which is invariant under $M$-action.
\end{theorem}

\vspace{2 mm}

Note that the original Bogolioubov's condition of amenability is contained in this theorem (see (3) and (4)).

\vspace{2 mm}

\underline{\textbf{Demonstration:}} Let us show (1). With its definition, it is clear that $M$ is a monoid. Let us define $D(f;g) = 
\int_K d(f(x);g(x))dP(x)$. $D$ is a pseudo-distance on $M$, and since $P$ is invariant under all maps of $M$, $D$ is invariant under 
right-translations. Then, if $(f_n)$ is a sequence of $M$, $(f_n)$ has a limit point in point-wise convergence topology (since $K^K$ 
is compact with Tychonov theorem), hence in $D$-topology since $(D(f_n;f)) \rightarrow 0$ thanks to Lebesgue theorem ($d$ is bounded 
on $K$ which is bounded too, hence the domination condition is satisfied). Now, we change our notations, and assume that $(f_n) 
\rightarrow f$ in point-wise convergence topology hence in $D$-topology. We have to show that $f \in M$. Let $g$ be a continuous map 
from $K$ to $\mathbf{R}$. One has $ \lim_{n\rightarrow + \infty}\int_K g\circ f_n (x) dP(x) = \int_K g\circ f (x) dP(x) = \int_K g(x) dP(x)$ with Lebesgue theorem ($g$ is continuous on a compact set hence bounded) and the invariance properties. Hence, since the 
continuous maps are dense in the set of integrable maps on compact sets, $\int_K g\circ f (x) dP(x) = \int_K g(x) dP(x)$ for any integrable 
map, and thus $P$ is invariant under $f$.

Now, let us show that $f \mapsto f(x;y)$ is measurable for almost every $x$ and $y$.
Denote $\mu$ a probability measure on $K$ such that $Supp(\mu) = K$, and : 
$$I(f;x;y) = \limsup_{r \rightarrow 0} \frac{\int_{B(x;r)}\int_{B(y;r)} d(f(t);f(s))d\mu(t)d\mu(s)}{\mu(B(x;r))\mu(B(y;r))}$$
for $(x;y) \in K$ and $f$ any Borel transformation of $K$.
$$d(f(t);f(s)) - d(g(t);g(s))\leq d(f(t);g(t)) + d(g(t);f(s)) - d(g(t);g(s))$$
$$\leq d(f(t);g(t)) + d(f(s);g(s))$$
for any $t$ and $s$. 
Hence, swapping $f$ and $g$,
$$\vert d(f(t);f(s)) - d(g(t);g(s))\vert\leq d(f(t);g(t)) + d(f(s);g(s)).$$
Thus, 
$$\frac{\vert\int_{B(x;r)}\int_{B(y;r)} d(f(t);f(s))d\mu^2(t;s)}{\mu(B(x;r))\mu(B(y;r))} - \frac{\int_{B(x;r)}\int_{B(y;r)} d(g(t);g(s))d\mu^2(t;s)\vert }{\mu(B(x;r))\mu(B(y;r))}$$
$$\leq \frac{2D(f;g)}{\mu(B(x;r))\mu(B(y;r))}$$
for any $r>0$. Hence $f \mapsto \frac{\int_{B(x;r)}\int_{B(y;r)} d(f(t);f(s))d\mu(t)d\mu(s)}{\mu(B(x;r))\mu(B(y;r))}$ is continuous for each $r>0$, and 
thus $f \mapsto I(f;x;y)$ is measurable for any $x$ and $y$. As a consequence, $(f;x;y) \mapsto I(f;x;y)$ is measurable

Now, since $f$ is measurable, $(x;y) \mapsto d(f(x);f(y))$ is measurable on $K^2$. Hence, with Lusin theorem, there is a null set $N$ such that $(x;y) \mapsto d(f(x);f(y))$ is continuous on $K^2 \setminus N$. Then : 
$$\vert I(f;x;y) - f(x;y) \vert \leq $$
$$\limsup_{r \rightarrow 0} \frac{\int_{B(x;r)\times B(y;r)} \vert d(f(t);f(s)) - d(f(x);f(y))\vert d\mu(t)d\mu(s) }{\mu(B(x;r))\mu(B(y;r))}$$.

Now, assuming that $(x;y) \in K^2 \setminus N$, $(t;s) \mapsto d(f(t);f(s))$ is continuous in $(x;y)$, and one easily sees that for $r$ small enough, $\vert I(f;x;y) - f(x;y) \vert $ is as small as desired. Hence $I(f;x;y) = d(f(x);f(y))$ for almost every $(x;y)$, and $(f;x;y) \mapsto f(x;y)$ is measurable.

The last assertion is a consequence of the Lebesgue theorem.

\vspace{2 mm}

Now, we show (2). Just take for $M$ the completion of the monoid generated by $S$ for the Cauchy sequences, and choose for each $[f]\in M\setminus S$ defined by some Cauchy sequence $(f_n) \in S$ a limit point $f\in K^K$ which is compact thanks to Tychonov theorem. Hence $(f_n) \rightarrow f$ in $D$-topology. Then $M$ satisfies the requested properties.

\vspace{2 mm}

(3) is just an application of our main theorem "Geometry and Probability".

\vspace{2 mm}

Now we show (4). Define $d'(x;y) = \int_M d(f(x);f(y))dm(f)$. We quickly see that $d'$ is a pseudo-distance. 
Let us show that $d'$ is invariant under $M$-action. 

$$d'(x;y) = \int_M d(f(x);f(y))dm(f) =$$
$$\int_M d(f(h(x));f(h(y)))dm(f) = d'(h(x);h(y))$$
for any $h\in M$ and every $(x;y)$, with invariance of $m$ under $M$-action.

Now, $(f;x;y) \mapsto d(f(x);f(y))$ is measurable since it is component-wisely. Hence, for $\mu$-almost every $(x;y)$, $f \mapsto d(f(x);f(y))$ is continuous almost everywhere. Hence, with Lebesgue theorem, for almost every $x$, 
$$\lim_{y \rightarrow x} \int_M d(f(x);f(y))dm(f) = 0.$$ 
Then, denote $Q$ a $d$-compact set such that for any $x\in Q$ the previous identity holds. Let $(x_n)$ be any sequence of $Q$. It admits a $d$-converging subsequence, hence a $d'$--converging subsequence. $Q$ is $d'$-compact. Thus, $K$ can be written as the union of a $\mu$-null set with a $\sigma$-compact set.

Let us show the last point : if $F \subset \Sigma$ is closed for $d'$, and if $F \ni (y_n) \rightarrow y\in K$ for $d$, we have already shown that $(d'(y_n;y))\rightarrow 0$. Hence $y \in F$ and $F$ is closed for $d$. Now, if $d'$ is a true distance, and if $F\subset \Sigma$ is closed for $d$, let $F \ni (y_n) \rightarrow y\in K$ for $d'$. $F$ is compact, hence $(y_n)$ has a subsequence $(z_n)$ converging to some $z \in F$. Then $(z_n) \rightarrow y = z\in F$ for $d'$ since $K$ is Hausdorff for $d'$, and $F$ is closed for $d'$.

Finally, let us show (5). Write $K = N \cup \bigcup_n Q_n$ where $N$ is $\mu$-null and $Q_n$ is $d'$-compact for each $n$. Use our main theorem to build a probability $P_n '$ on each $Q_n$ such that it is invariant under partial $d'$-isometries. Prolonging $P_n '$ by zero, we define new probabilities $P_n$ on $K$, and these probabilities are invariant under $S$-action when taken for Borel sets such that they and their image are contained in their support. Now $(P_n)$ is a sequence of $d'$-Borel probabilities defined on a compact set, hence it has a limit point $P$ in weak topology. Then $P$ is invariant under $S$-action. Indeed, for any $f \in S$, any Borel set can be written as the denumerable union of compact subsets of $Q_n \cap f^{-1}(Q_n)$ with a $P$-null set $N$ such that $f(N)$ is null.\begin{flushright}
$\blacksquare$
\end{flushright}

\section*{Conclusion}

In this paper, we have seen that general theorems dealing with the question of existence of probabilities invariant by a groupoid of maps exist. Indeed, we defined three different notions of compatibility of a measure with a metric structure, and we showed the existence of such probabilities on any compact metric space. Some results about their uniqueness are reachable. The topological nature of the underlying space is the key point for this study. More generally, the properties of the spaces of compatibles probabilities on a compact metric space contain informations about its metric structure, peculiarly about its homogeneity. 

Our construction of such compatible probabilities needs the choice of a ultra-filter of neighbourhoods of $0$ in $(0;+ \infty )$, and consequently the axiom of choice has to be assumed. Once a ultra-filter is chosen, some functoriality results can be shown in the category of compact metric spaces endowed with the probabilities we built.

\vspace{2 mm}

This work uses a strategy for the building of probability measures which is very close to the original method of Haar for his own probability. Nevertheless, we managed to show really new results generalising in some particular cases the well-known consequences of the Haar work. The point is that in metric cases, the group assumption is non longer needed.

\vspace{2 mm} 
 
We illustrated our theorem with many consequences and applications. First of all, we exhibited the link of our construction with the usually known measures : the Lebesgue one, the Haar one and the Hausdorff one. We showed some generalisations of the "Geometry and Probability Theorem" too. And then, we gave some direct applications, notably in a concrete case in Phylogeny.

\vspace{2 mm}

Finally, our most impressive consequence is probably the Krylov - Bogolioubov Equivalence Theorem, which gives very precise informations for the existence or not of a probability invariant by a monoid of measurable maps acting on a compact metric space.
 
\vspace{2 mm} 

New perspectives appear to follow our paper. First of all, a more general theorem must be reachable in the framework of uniform compact Hausdorff spaces. For example, one can try to show such a theorem in compact spaces $K$ endowed with some kind of distance satisfying the following axiomatic: 

\vspace{2 mm}

There is a totally ordered abelian group $D$ and a map $ d : K^2 \rightarrow D$ such that:

- $d$ is non-negative,

- $d(x;y) = 0_D \Leftrightarrow  x=y$,

- $d(x;y) = d(y;x)$,

- $d(x;y) \leq d(x;z) + d(z; y)$,

- the topology of $K$ is generated by the $d$-balls.

\vspace{2 mm}

We are nearly convinced that almost all uniform compact Hausdorff spaces can be endowed with such a "distance". And this structure defines the groupoid of partial "isometries" which is needed to formulate the original problem of Integral Geometry.

\vspace{2 mm}

Such a theorem would generalise the most general case of the Haar measure (when the Hausdorff compact group is not metrizable but uniform). This theorem would then be another good answer for the Integral Geometry problem.

\vspace{2 mm}

Note that the previous approach enables to state some interesting theorems for Functional Analysis just using the well known Haar measure.

Indeed, take for $E$ the space of continuous functions on a compact metric space $(K;d)$ to $\mathbf{R}$, endowed with the following generalized distance.

\vspace{2 mm}

Let $f \in E$. $f$ is uniformly continuous, hence the map $$\omega (f)(\varepsilon ) = \sup_{(x;y) \in K^2, d(x;y)\leq \varepsilon } \vert f(x)-f(y)\vert$$ is defined and decreasing, and satisfies : $\lim_{\varepsilon \rightarrow 0} \omega (f)(\varepsilon ) =0$. Conversely, if $\lim_{\varepsilon \rightarrow 0} \omega (f)(\varepsilon ) =0$, $f$ is continuous.

Now, the set of decreasing maps $g : (0;+\infty) \rightarrow \mathbf{R}$ such that $\lim_{\varepsilon \rightarrow 0} g(\varepsilon ) =0$ can be mapped into a totally ordered abelian group. Indeed, choose a ultra-filter $U$ containing the filter of neighbourhoods of $0$ in $(0;+\infty)$, and consider the set $D$ of equivalence classes of maps $g$ for the equality along some elements of $U$. Showing that $D$ is a totally ordered abelian group is straightforward. Denote $[g]$ the equivalence class of $g : (0;+\infty) \rightarrow \mathbf{R}$ for that relation.

Then, we set the following definition for a generalized distance on $E$ :

$$\delta (f;g) = [\sup_{x \in K} \vert f(x) - g(x)\vert\times Id_{(0;+\infty)}] + [\omega (f-g)].$$

Checking that it is a generalized distance raises no problem.

Now, endow $E$ with the topology generated by $\delta$, and note that $E$ is a topological group for the addition law.

In fact, $E$ is locally compact. Indeed, for any $f \in E$ and any $\Delta \in D$, $\overline{B(f;\Delta)}$ is a closed bounded subset of $E$, and it is equicontinuous. Hence with Ascoli's theorem, it is compact.

Consequently, $E$ admits a Haar measure invariant by all translations. One can choose that Haar measure such that the measure of $\overline{B(0;1)}$ equals $1$.

\vspace{2 mm}

The same kind of constructions must give results looking like ours when no group structure is assumed.

\vspace{4 mm}

Last point, but not least, the formula of theorem $4$ gives an explicit definition for the measure $\mu$. Generalising that formula, one can define some new notions of probabilities. For example, one can define a "generalized probability" by the following axioms:

- a space $E$ and a $\sigma$-algebra $B$ on $E$ are given,

- a commutative ordered monoid $(F;+)$ is given, in which suprema and infima always exist, and whose neutral element is the lowest,

- a map $P : B \rightarrow F$ is given, called "generalized probability", with:

	- $P(E) > 0_F$,
	
	- If $(A_n)$ is a countable family of disjoint sets of $B$, 
	$$P(\bigcup_n A_n) = \sup_{finite families} \sum_k P(A_k) .$$
	
\vspace{2 mm}	

Assume then that such a monoid $F$ and an increasing morphism from the set of bounded sequences of rational numbers to $F$ are given, the image of a sequence $(r_n)$ being independent from the first values of $(r_n)$. Denote such a morphism $L$. Then, in a compact metric space, the following version of the formula of theorem $4$:

\begin{equation*}
\mu_K(A) = \sup_{ F \subseteq A, F closed} \inf_{\delta >
0} L(( \frac{N(\overline{F}^{\delta};\frac{1}{n})}{%
N(K;\frac{1}{n})})_n) 
\end{equation*}

gives a candidate for a generalised probability in $F$. Non trivial examples for such couples $(F;L)$ exist: one can take, as in this paper, $F = \mathbf{R}$ and $L$ a limit notion. Another example is to take for $F$ the Dedekind completion of ${}^* \mathbf{R}$, the field of non standard real numbers, together with the morphism $L$ defined on the ring of bounded sequences of real numbers and mapping into ${}^* \mathbf{R}$. Unfortunately, this last example does not seem very interesting since the Dedekind completion of ${}^* \mathbf{R}$ does not have pleasant properties. 

Further demonstrations are needed, and a better approach may exist, but these ideas suggest that generalisations of the actual Probabilities Theory are reachable, and that constructions of such generalized probabilities are possible. Such generalised theories could be good frameworks to define some new kinds of "measures" with good invariance properties in difficult cases of Functional Analysis dealing with separable non-locally compact spaces.

\vspace{15 mm}

{\Large Greeting} I thank M. Paulin (Universit\'{e} Paris-Sud), M. Penot
(Universit\'{e} des pays de l'Adour) and M. Ciligot-Travain (Universit\'{e}
d'Avignon) for the interest they manifested for my work. I express my
gratitude to Mme Kosmann-Schwarzbach, M. Aubin (Professeurs honoraires),
M. Zvonkine (Universit\'{e} de Bordeaux I) and specially to M. Lewandowski (Université de Rennes) and M. Léonard (Université Paris-Ouest) for their encouragements and their help. Finally, I
warmly thank M. Drouin and M. Paya, Teacher like me at Lyc\'{e}e de Borda in Dax, for
their work of re-lecture and check of demonstrations, and for our fruitful discussions.


\clearpage

\begin{thebibliography}{99}

\bibitem{calculproba} Foata Dominique, Fuchs Aim\'{e} \emph{Calcul des
probabilit\'{e}s}. Dunod, deuxi\`{e}me \'{e}dition, 1998.

\bibitem{topobou} Bourbaki Nicolas \emph{Topologie g\'{e}n\'{e}rale}.
Springer, r\'{e}impression inchang\'{e}e de l'\'{e}dition originale.

\bibitem{Bourbaki mesure} Bourbaki Nicolas \emph{Int\'{e}gration Chapitre 9}%
. Springer, r\'{e}impression inchang\'{e}e de l'\'{e}dition originale.

\bibitem{topovect} Bourbaki Nicolas \emph{Espaces vectoriels topologiques}.
Springer, r\'{e}impression inchang\'{e}e de l'\'{e}dition originale.

\bibitem{Fractales} Tricot Claude \emph{G\'{e}om\'{e}tries et mesures
fractales : une introduction}. Ellipses, 2008.

\bibitem{Brézis} Br\'{e}zis Ha\"{\i}m \emph{Analyse fonctionnelle}. Dunod,
1999.

\bibitem{Rudin} Rudin Walter \emph{Analyse r\'{e}elle et complexe}. Dunod, Sciences Sup,1998.

\bibitem{Laudenbach} Laudenbach Fran\c{c}ois \emph{Calcul diff\'{e}rentiel
et int\'{e}gral}. \'{E}ditions de l'\'{E}cole Polytechnique, 2000.

\bibitem{vill1} Otto F. and Villani C. \emph{Generalization of an inequality by Talagrand and links with the logarithmic Sobolev inequality}. J. Funct. Anal. 173, 2 (2000), 361-400.

\bibitem{vill2} Villani C. \emph{Optimal transport, old and new. Notes de cours pour l'Ecole d'été de Saint-Flour 2005,} see : \underline{www.umpa.ens-lyon.fr/~villani}.

\bibitem{tree} S. N. Evans \emph{Probability and real trees. Ecole d'été de Probabilités de Saint-Flour XXXV - 2005}. Springer


\vspace{3 mm}

{\large Ultrafilters and non standard analysis}

\vspace{3 mm}

\bibitem{cech} Cech E. \emph{On bicompact spaces}, Ann. Math. 38: 823-844, 1937.

\bibitem{stone} Stone M. H. \emph{Applications of the theory of Boolean rings to general topology},  Trans. Amer. Soc. 41(3): 375-481, 1937.

\bibitem{robinson} Robinson A \emph{Non-standard analysis}, North-Holland Publishing Co., Amsterdam 1966.

\vspace{3 mm}

{\large Invariant measures and generalisations}

\vspace{3 mm}

\bibitem{haar} Haar A. \emph{Der Massbegriff in der theorie der kontinuierlichen Gruppen} Annals of Mathematics, 2 34(1) : 147-169, 1933.

\bibitem{cartan} Cartan H. \emph{Sur la mesure de Haar} Comptes-rendus de l'Académie des Sciences, vol. 211:759-762, 1940.

\bibitem{hausdorff} Hausdorff F. \emph{dimension und \"ausseres Mass} Mathematische Annalen 79 (1-2) : 157-179, 1918.

\bibitem{falconer} Falconer K. J. \emph{The geometry of fractal sets} Cambridge University Press, 1985.

\bibitem{contrex} Taylor \emph{The exact Hausdorff measure of the sample path for planar Brownian motion}. Proceedings of the Cambridge Philosophical Society, 1964.

\bibitem{loom} Loomis L. H. \emph{Haar measure in uniform structures}, Duke Math. J. 16: 193-208, 1949.

\bibitem{Chris} Christensen J. P. R. \emph{On some measures analogous
to Haar measures}. Math. Scand. 26, 103-106 1970.

\bibitem{Brownell} Brownell F. H. \emph{Translation invariant measure over separable Hilbert spaces and other translation spaces}, Pacific J. Math, 3: 531-553, 1953. 

\bibitem{cartier-morette} Cartier P., De Witt-Morette C. \emph{A new perspective on Functional Integration}, Journal of Math. Phys., 36:2137-2340, 1995.

\bibitem{camo} Cartier P., De Witt-Morette C. \emph{Intégration fonctionnelle; éléments
d'axiomatique} C.R. Acad. Sci. Paris, t.316, Série II, 733-738, 1993.

\bibitem{gross} Gross L. \emph{Abstract Wiener space} in: Proc. 5th Berkeley Symp. Math. Stat. and Probab. 2, part 1: 31-42, University of California Press, Berkeley, 1965.

\bibitem{wien} Wiener N. \emph{Differential space} J. Math. and Phys., 2: 132-174, 1923.

\bibitem{Jia} Yan J-A \emph{Generalizations of Gross's and Minlos's theorems}, Springer-Verlag, in: Séminaire de probabilités (Strasbourg), tome 23, 395-404, 1989. 

\bibitem{sudakov} Sudakov, V. N., \emph{Linear sets with quasi-invariant measure.} (in Russian), Doklady Akademii Nauk 127, 524-525, 1959.

\bibitem{krybog} N. N. Bogoliubov, N. M. Krylov, \emph{La theorie generale de la mesure dans son application à l'étude de systèmes dynamiques de la mecanique non-lineaire} (in French), Ann. Math. II, 38 (1): 65?113, 1937.


\end{thebibliography}
\end{document}